\documentclass[twoside,draft]{article}

\usepackage{amsfonts}
\usepackage{stmaryrd}
\usepackage{amssymb}
\usepackage{euscript}
\usepackage{amsthm}
\usepackage{amsmath}
\usepackage{amscd}
\usepackage{latexsym}
\usepackage{mathrsfs}
\usepackage{graphicx}
\usepackage{color}
\usepackage{dsfont}
\usepackage{bm}
 \usepackage{cleveref}
\usepackage{enumitem}%

\numberwithin{equation}{section}

\newtheorem{theorem}{Theorem}[section]
\newtheorem{lem}{Lemma}[section]
\newtheorem{pro}{Proposition}[section]
\newtheorem{cor}{Corollary}[section]
\newtheorem{rem}{Remark}[section]
\newtheorem{rems}{Remarks}[section]
\newtheorem{ex}{Example}[section]
\newtheorem{defi}{Definition}[section]
\newtheorem{hyp}{Assumption}[section]

\newcommand{\bt}{\begin{theorem}}
\newcommand{\et}{\end{theorem}}
\newcommand{\bl}{\begin{lem}}
\newcommand{\el}{\end{lem}}
\newcommand{\bp}{\begin{pro}}
\newcommand{\ep}{\end{pro}}
\newcommand{\bcor}{\begin{cor}}
\newcommand{\ecor}{\end{cor}}

\newcommand{\bd}{\begin{defi} \rm }
\newcommand{\ed}{\end{defi}}
\newcommand{\brem }{\begin{rem} \rm }
\newcommand{\erem }{\end{rem}}
\newcommand{\brems }{\begin{rems} \rm }
\newcommand{\erems }{\end{rems}}
\newcommand{\bhyp }{\begin{hyp} \rm }
\newcommand{\ehyp }{\end{hyp}}
\newcommand{\bex}{\begin{ex} \rm }
\newcommand{\eex}{\end{ex}}

\newcommand{\newA}{A}
\newcommand{\newK}{K}

\newcommand{\gen}{g}
\newcommand{\Gen}{G}

\newcommand{\wWl}{W^{\lambda}}
\newcommand{\xXl}{X^{\lambda}}

\newcommand{\ee}{e}
\newcommand{\yxi}{\xi}
\newcommand{\bxi}{\bar{\xi}}

\newcommand{\seq}[1]{{\lbrace #1 \rbrace}}

\newcommand{\bigcdot}{\mathbin{{\hbox{\scalebox{.65}{$\bullet$}}}}}
\newcommand{\bcd}{\bigcdot}

\newcommand{\wt}{\widetilde}
\newcommand{\wh}{\widehat}
\newcommand{\ov}{\overline}

\newcommand{\I}{\mathds{1}}

\newcommand{\llb}{\llbracket}
\newcommand{\rrb}{\rrbracket}

\newcommand{\esssup}{\operatornamewithlimits{ess\,sup}}
\newcommand{\essinf}{\operatornamewithlimits{ess\,inf}}

\newcommand{\xtau}{\sigma}
\newcommand{\xnu}{\nu}

\newcommand{\STF}{{\cT}(\FF)}
\newcommand{\STFt}{{\cT}_{t,T}(\FF)}
\newcommand{\STFtau}{{\cT}_{\tau,T}(\FF)}
\newcommand{\STFnu}{{\cT}_{\nu,T}(\FF)}
\newcommand{\STFp}{{\cT}^p(\FF)}

\newcommand{\STFtaup}{{\cT}^+_{\tau,T}(\FF)}

\newcommand{\bSTF}{\overline{\cT}(\FF)}
\newcommand{\bSTFt}{\overline{\cT}_{t,T}(\FF)}
\newcommand{\bSTFp}{\overline{\cT}^p(\FF)}

\newcommand{\bSTFtau}{\overline{\cT}_{\tau,T}(\FF)}

\newcommand{\bSTFtaup}{\overline{\cT}^+_{\tau,T}(\FF)}

\newcommand{\yyy}{y}
\newcommand{\zzz}{z}
\newcommand{\kkk}{k}

\newcommand{\ww}{w}

\newcommand{\yy}{y}

\newcommand{\xS}{S}

\newcommand{\cKS}{\cK(\xS)}
\newcommand{\bcKS}{\ov{\cK}(\xS)}
\newcommand{\bcK}{\ov{\cK}}

\newcommand{\bzeta}{\bar{\zeta}}

\newcommand{\barR}{\ov{R}}
\newcommand{\barS}{\ov{S}}

\newcommand{\barV}{\ov{V}}
\newcommand{\barY}{\ov{Y}}
\newcommand{\barZ}{\ov{Z}}
\newcommand{\barK}{\ov{K}}

\newcommand{\whY}{\wh{Y}}

\newcommand{\whZ}{\wh{Z}}

\newcommand{\whK}{\wh{K}}

\newcommand{\wtzeta}{\wt{\zeta}}
\newcommand{\wtxi}{\wt{\xi}}
\newcommand{\wtY}{\wt{Y}}
\newcommand{\wtZ}{\wt{Z}}
\newcommand{\wtK}{\wt{K}}

\newcommand{\wtU}{\wt{U}}
\newcommand{\tPP}{\widetilde{\mathbb P}}

\newcommand{\wtGamma}{\widetilde{\Gamma}}

\newcommand{\cH}{\mathcal{H}}
\newcommand{\cS}{\mathcal{S}}
\newcommand{\cK}{\mathcal{K}}

\newcommand{\cOff}{\mathcal{O}(\FF)}

\newcommand{\cPff}{\mathcal{P}(\FF)}

\newcommand{\cE}{{\mathcal E}}

\newcommand{\cF}{{\mathcal F}}
\newcommand{\cG}{{\mathcal G}}

\newcommand{\cT}{{\mathcal T}}

\newcommand{\FF}{{\mathbb F}}
\newcommand{\GG}{{\mathbb G}}
\newcommand{\RR}{{\mathbb R}}

\newcommand{\PP}{{\mathbb P}}
\newcommand{\NN}{{\mathbb N}}
\newcommand{\EE}{{\mathbb E}}
\newcommand{\EP}{\mathbb{E}_{\mathbb{P}}}

\title{{\Large \bf  WELL-POSEDNESS AND PENALIZATION SCHEMES FOR GENERALIZED BSDEs AND REFLECTED GENERALIZED BSDEs} \vskip 35 pt }
\author{Libo Li$\,^{a}$, Ruyi Liu$\,^{b}$ and Marek Rutkowski$\,^{b,c}$ \\ \\ \\ \\
\\ $^{a\,}$School of Mathematics and Statistics, University of New South Wales \\ Sydney, NSW 2052, Australia \\ \\
$^{b\,}$School of Mathematics and Statistics, University of Sydney \\ Sydney, NSW 2006, Australia \\ \\
$^{c\,}$Faculty of Mathematics and Information Science, Warsaw University of Technology \\ 00-661 Warszawa, Poland \\ }

\date{\vskip 45 pt \today \vskip 35 pt}

\begin{document}

\maketitle

\begin{abstract}
The paper is directly motivated by the pricing of vulnerable European and American options in a general hazard process setup and a related study of the corresponding pre-default backward stochastic differential equations (BSDE) and pre-default reflected backward stochastic differential equations (RBSDE). We work with a generic filtration $\FF$ for which the martingale representation property is assumed to hold with respect to a square-integrable martingale $M$ and the goal of this work is of twofold. First, we aim to establish the well-posedness results and comparison theorems for a generalized BSDE and a reflected generalized BSDE with a continuous and nondecreasing driver $A$. Second, we study extended penalization schemes for a generalized BSDE and a reflected generalized BSDE in which we penalize against the driver in order to obtain in the limit either a particular optimal stopping problem or a Dynkin game in which the set of admissible exercise time is constrained to the right support of the measure generated by $A$.

\end{abstract}

\newpage
\tableofcontents
\newpage

%%%%%%%%%%%%%%%%%%%%%%%%%%%%%%%%%%%%%%%%%%%%%%%%%%%%%%%%%%%%%%%%%%%%%%%%%%%%%
%%%%%%%%%%%%%%%%%%%%%%%%%%%%%%%%%%%%%%%%%%%%%%%%%%%%%%%%%%%%%%%%%%%%%%%%%%%%%
\section{Introduction}    \label{sec1}
%%%%%%%%%%%%%%%%%%%%%%%%%%%%%%%%%%%%%%%%%%%%%%%%%%%%%%%%%%%%%%%%%%%%%%%%%%%%%
%%%%%%%%%%%%%%%%%%%%%%%%%%%%%%%%%%%%%%%%%%%%%%%%%%%%%%%%%%%%%%%%%%%%%%%%%%%%%

In this work, we consider the generalized backward stochastic differential equation (GBSDE) and the reflected generalized backward stochastic differential equation (RGBSDE). More specifically, given a continuous, nondecreasing process $A$, we are interested in a generalized BSDE of the form
\begin{align} \label{sec1eq1}
Y_t = \xi_T - \int_{]t,T]} Z_s\,dM_s + \int_{]t,T]} g(s,Y_s)\,dA_s
\end{align}
where $M$ is a square-integrable martingale with respect to a reference filtration $\FF$ and $A$ is an adapted, continuous, nondecreasing process.
Furthermore, we study a reflected generalized BSDE given by
\begin{align} \label{sec1eq2}
Y_t = \xi_T - \int_{]t,T]} Z_s\,dM_s + \int_{]t,T]} g(s,Y_s)\,dA_s + K_T- K_t
\end{align}
where $\zeta$ is the lower obstacle and the nondecreasing process $K$ satisfies the appropriate Skorokhod conditions (see \eqref{RGBSDE}). To avoid confusion, we shall refer to the process $A$ as the {\it driver} whereas, as usual, the mapping $g$ is called the {\it generator} of a GBSDE.

Our motivation for studying equations \eqref{sec1eq1} and \eqref{sec1eq2} comes from the study of a BSDE and RBSDE in the progressive enlargement
$\GG$ of $\FF$ through observations of a random time, which can be interpreted as a default time of some credit-risky entity, and their subsequent applications to the pricing of vulnerable options of either a European or an American style. In particular, the form of the driver $A$ in equations  \eqref{sec1eq1} and \eqref{sec1eq2} is not arbitrary since it is directly related to the choice of a particular model for the default time, as was studied, e.g., in a recent paper by Jeanblanc and Li \cite{JL2020} and numerous previous works on credit risk modeling.  A study of a BSDE in the progressively enlarged filtration $\GG$ has been explored in two main directions. On the one hand, one can work directly in the enlarged filtration $\GG$, as done in Blanchet-Scalliet et al. \cite{BER2010}, Eyraud-Loisel and Royer-Carenzi \cite{ER2010}, Grigorova et al. \cite{GQS2020} and Dumitrescu et al. \cite{DQS2018}. On the other hand, in papers by Kharroubi and Lim \cite{KL2014} and Cr\'epey and Song \cite{CS2016}, the authors developed an alternative approach, which hinges on reducing the BSDE in the filtration $\GG$ to a reference filtration $\FF$ and subsequently examining the {\it $\FF$-reduced} BSDE, which is also known as the {\it pre-default} BSDE in applications of the theory of BSDEs to credit risk models.

The present work is motivated by the above-mentioned $\FF$-reduction approach and stems from the recent works of Aksamit et al. \cite{ALR2021} and Li et al. \cite{LLR1} where the pre-default GBSDE as well as the pre-default RGBSDE were examined and applied to solve super-hedging problems for vulnerable European and American options within a fairly general setup. In particular, the driver process $A$ appearing in the pre-default GBSDE and RGBSDE satisfied by the pre-default price of a vulnerable European or American option corresponds to the hazard process associated with the default time while the generator $g$ is used to model the uncertainty in the hazard process.

%Although we have assumed that the process $A$ is continuous it is not difficult to extend our results to the case where $A$ is allowed to have ordered jumps through a backward jump adapted methodology.

%with respect to the indicator $\I_{\llb \nu, \infty \llb}$ where $\nu$ is a random time which is not a stopping time in the reference filtration $\FF$.

%In the literature of default risk, the filtration enlargement is a fundamental approach to study the issue of defaultable claims. Introducing the default times  $\nu$ which are not $\FF$-stopping times into the model, the reference filtration $\FF$ would be progressively enlarged, by observing the occurrences of $\nu$, to its smallest extended filtration $\GG$ which accommodates $\nu$ as its stopping times. This approach was first introduced in fundamental works by Jeulin \cite{J1980}, Jeulin and Yor \cite{JY1978} and Yor \cite{Y1978}. Afterwards, many works, such as Aksamit et al. \cite{ACJ2015},  Choulli et al. \cite{CDV2020}, Duffine and Singleton \cite{DS2003} and many others, have worked on this subject especially in applications.

%differences with Aksamit et al.\cite{ALR2021}
%As a supplement to the well-posedness of the reduction for $\GG$ BSDE and  $\GG$ reflected BSDE in Aksamit et al.\cite{ALR2021}, we prove the uniqueness and existence in cases of the continuous and square-integrable hazard process $\Gamma$.

% Such a reduced $\FF$ BSDE was driven by a l\`adl\`ag process and the existence result was obtained by using a particular family of BSDE.

The existing literature on (generalized) BSDE and RBSDE are numerous and a non-exhaustive list include Pardoux and Zhang \cite{PZ1998}, Ren and El Otmani \cite{RE2010}, Essaky and Hassani \cite{EH20111, EH20112}, Essaky et al. \cite{EHO2013, EHO2015} and Eddahbi et al. \cite{EFO2017}. We mention also Grigorova et al. \cite{GIOOQ2017, GIOQ2018, GIOQ2020}, Hamad\`ene et al. \cite{HO2016, HHO2010} and Klimsiak et al. \cite{KRS2019}. However, to the best of our knowledge, no single piece of work encompasses the general framework which is adopted in papers by Aksamit et al. \cite{ALR2021} and Li et al. \cite{LLR1} and hence we have found it necessary to study the generalized BSDE within the framework consistent with the setup of \cite{LLR1}.

Our first goal is to obtain, under Assumptions \ref{nass1.1} and \ref{nass1.2} complemented by the postulate that the driver $A$ is bounded, the well-posedness and comparison results for a GBSDE and RGBSDE, as given by \eqref{sec1eq1} and \eqref{sec1eq2}, respectively. In particular, we establish
comparison results in Propositions \ref{npro1.1} and \ref{npro1.5}, while the existence and uniqueness results are given in Propositions \ref{npro1.3} and \ref{npro1.6}. Furthermore, we show in Proposition \ref{npro1.4} that for a GBSDE the boundedness assumption on $A$ can be lifted and
the existence and uniqueness result can be extended to the case where $A$ is square-integrable under the additional assumption that the generator $g(t,y)$ is non-negative and nondecreasing in the variable $y$. We expect the same technique to be applicable to an RGBSDE and hence omit the details here.

Our second goal is to present the penalization scheme for \eqref{sec1eq1} and \eqref{sec1eq2} where we penalize against the driver $A$. Similar to the classical case where $A_t = t$, the limiting process obtained from the penalization scheme correspond to either the value process of an optimal stopping problem (see Theorems \ref{th3.1} and \ref{th6.1}) or the value of a two player zero-sum game (see Theorem \ref{th6.2}). Compared to existing results on penalization of an RBSDEs, such as Hamad\`ene et al. \cite{HHO2010}, the new feature here is that the set of admissible exercise times is restricted to the right support of the random measure generated by the driver $A$.

 Finally, it is worth mentioning that the assumption of continuity of $A$ is not essential and thus, at least in principle, in can be relaxed. In particular, it would be possible to work also with a purely discontinuous driver by considering a left-continuous nonlinear driver and hence a l\`agl\`ad GBSDE, as opposed to a c\`adl\`ag GBSDE considered in the present work. Then the set of admissible exercise times for the optimal stopping problem, which can be obtained by considering the limiting process of the penalization scheme, would be restricted to the jump times of $A$ (in financial applications, this would formally correspond to the case of Bermuda options). An extension of the setup and results obtained in the present work to the case of a general driver is the topic of our ongoing research.

\section{Preliminaries}         \label{nsec1}
%%%%%%%%%%%%%%%%%%%%%%%%%%%%%%%%%%%%%%%%%%%%%%%%%%%%%%%%%%%%%%%%%%%%%%%%%%%%%%%%%%%%%%%%%%%

Our goal is to analyze some classes of generalized BSDEs and reflected generalized BSDEs with respect to the filtration $\FF$.
We work throughout under the following set of standing assumptions, which specify the inputs in a generalized BSDE (GBSDE).

\bhyp \label{nass1.1}
We assume that we are given the following objects:\\
(i) a probability space $(\Omega,{\cal F},\PP )$ endowed with a filtration $\FF$ satisfying the usual conditions of
  right-continuity and completeness; \\
(ii) an $\RR^{d}$-valued, square-integrable martingale $M$, which is assumed to have the predictable representation property for the filtration $\FF$; \\
(iii) an $\FF$-adapted, nondecreasing and continuous process $\newA$, which is square-integrable so that $\EE[\newA_T^2] < \infty$.
By convention, we set $A_{0-}=A_0=0$; \\
(iv) an ${\cal F}_T$-measurable and square-integrable random variable $\xi_T$.
\ehyp

We introduce the space $\cH^2(M)=\{ Z \in \cPff : \|Z\|^2_{\cH^2(M)}< \infty \}$ with the pseudo-norm
\begin{align*}
\|Z\|_{\cH^2(M)}:=\Big[ \EE\big[(Z^2\bigcdot [M])_T\big] \Big]^{1/2}
\end{align*}
where the martingale $M$ is assumed to be square-integrable so that $\EE [[M]_T] < \infty$.
Let $\STF$ denote the class of all $\FF$-stopping times with values in $[0,T]$. It is known that the space
$\cS^2(\FF)=\{ X\in\cOff :\|X\|^2_{\cS^2(\FF)}< \infty \}$ with the norm
\begin{align*}
\|X\|_{\cS^2(\FF)}:=\bigg[ \EE\bigg(\esssup_{\tau \in \STF}|X_\tau|^2\bigg)\bigg]^{1/2}
\end{align*}
is a Banach space (see Proposition 2.1 in Grigorova et al. \cite{GIOOQ2017}).

We denote by  $\cK$ (resp., $\bcK$) the class of all c\`adl\`ag, nondecreasing, $\FF$-predictable (resp., l\`adl\`ag, nondecreasing, $\FF$-predictable) processes. Recall that a stochastic process $X$ with sample paths possessing right-hand limits is said to be {\it $\FF$-strongly predictable} if $X$ is $\FF$-predictable and the process $X_+$ is $\FF$-optional (see Definition 1.1 in  Gal'\v{c}uk \cite{G1982}).
In particular, all processes from the class $\bcK$ is $\FF$-strongly predictable.
% From Theorem 1.33 in Jacod \cite{Jxxx} (or Dellacherie IV T41 \cite{Dxxx})
Furthermore, any process $K \in \cK$ with $K_0=0$ has a unique decomposition $K=K^c+K^d$ where $K^c_0=K^d_0=0$, $K^c$ is an $\FF$-adapted, continuous, nondecreasing process and $K^d$ is an $\FF$-predictable, c\`adl\`ag, purely discontinuous, nondecreasing process. More generally, if $K$ belongs to $\bcK$ and $K_0=0$ then the decomposition becomes $K=K^c+K^d+K^g$ where $K^g$ with $K^g_0=0$ is an $\FF$-adapted, c\`agl\`ad, purely discontinuous, nondecreasing process. If $X$ and $Y$ are arbitrary $\FF$-optional processes, then the inequality $Y\geq X$ means that $Y_{\xtau}\geq X_{\xtau}$ for every $\xtau \in \STF$.

For future reference, we recall the classical version of the Doob-Meyer-Mertens decomposition theorem (see Mertens \cite{M1972} and El Karoui \cite{EK1981}) and the basic properties of the Snell envelope. We stress that the filtration $\FF$ is assumed to satisfy the usual conditions of Assumption \ref{nass1.1} (i). The interested reader is referred to Gal'\v{c}uk \cite{G1982} for an extended version of the Doob-Meyer-Mertens decomposition where these assumptions about the filtration $\FF$ are relaxed.

Let us recall the notion of a \emph{strong supermartingale,} which is known to coincide with the classical concept of a \emph{supermartingale} if a process is assumed to be c\`adl\`ag. We will also use later the concept of a \emph{strong $\cE^g$-supermartingale}
where $\cE^g$ denotes the nonlinear evaluation generated by solutions to a (generalized) BSDE with generator $g$ (see Peng \cite{P2004a,P2004}).

\bd \label{ndef1.1}
A process $X$ is called a \emph{strong supermartingale} if it is $\FF$-optional and for all $\tau,\nu \in \STF$ such that $\tau\le\nu$ we have
\begin{equation} \label{neq3b}
X_\tau \ge \EE[ X_\nu \,|\, \cF_\tau ].
\end{equation}
\ed

\brem \label{nrem1.1}
It is well known that any strong supermartingale is a l\`adl\`ag process, that is, it has almost all sample paths with right-hand and left-hand limits so that the processes $X_{-}$ and $X_{+}$ are well-defined.  Moreover, any strong $\FF$-supermartingale is a right-upper-semicontinuous process so that $X \ge X_+$.
\erem

\bd  \label{ndef1.2}
We say that $X$ is the \emph{Snell envelope} of an $\FF$-optional process $\xi$ if: \\
(i) $X$ is a strong supermartingale such that $X\geq\xi$; \\
(ii) for any strong supermartingale such that $Y\geq\xi$ the inequality $Y\geq X$ holds.
\ed

\brem \label{nrem1.2}
The following simple observations will be useful: \\
(i) if $X$ is a strong supermartingale, then the Snell envelope of $X$ is equal to $X$; \\
(ii) the $\FF$-Snell envelope is monotone, in the sense that if $\eta \geq \xi$, then the Snell envelope of $\eta$ dominates the Snell envelope of $\xi$.
\erem

\bt \label{nth1.1}
Any strong supermartingale $Y$ of class (D) has the unique  Doob-Meyer-Mertens decomposition $Y=Y_0+N-B^c-B^d-B^g$ where $N$ is a uniformly integrable martingale, $B^c$ is an $\FF$-adapted, nondecreasing, continuous process, $B^d$ is an $\FF$-predictable, nondecreasing, purely discontinuous process, and $B^g$ is an $\FF$-adapted, c\`agl\`ad, nondecreasing, purely discontinuous process. Furthermore, $N_0=B^c_0=B^d_0=B^g_0=0$.
\et

%
%\brem \label{nrem4.1}
%Any l\`adl\`ag, $\FF$-predictable, nondecreasing process $B$ has a unique decomposition $B=B^c+B^d+B^g$ with the properties stated in Proposition \ref{pro4.1}. Observe that if we denote $A:=B^c+B^d$ and $C:=B^g_+$, then from  Proposition  \ref{pro4.1} we obtain the DMM decomposition $Y=N-A-C_{-}$,
%which is the form used in some papers on optimal stopping, e.g., \cite{EK1981,GIOOQ2017,GIOQ2020,KQ2012}. Alternatively, if we set $D=A+C_{-}$, then we obtain the decomposition $Y=N-D$ where $D$ is a nondecreasing, $\FF$-{\it strongly predictable} process, in the sense that $D$ is $\FF$-predictable and $D_{+}$ is $\FF$-optional (see \cite{G1981,G1982}).
%\erem

%%%%%%%%%%%%%%%%%%%%%%%%%%%%%%%%%%%%%%%%%%%%%%%%%%%%%%%%%%%%%%%%%%%%%%%%%%%%%%%%%%%%%%%%%%
%%%%%%%%%%%%%%%%%%%%%%%%%%%%%%%%%%%%%%%%%%%%%%%%%%%%%%%%%%%%%%%%%%%%%%%%%%%%%%%%%%%%%%%%%%
\section{Generalized BSDEs}         \label{nsec2}
%%%%%%%%%%%%%%%%%%%%%%%%%%%%%%%%%%%%%%%%%%%%%%%%%%%%%%%%%%%%%%%%%%%%%%%%%%%%%%%%%%%%%%%%%%
%%%%%%%%%%%%%%%%%%%%%%%%%%%%%%%%%%%%%%%%%%%%%%%%%%%%%%%%%%%%%%%%%%%%%%%%%%%%%%%%%%%%%%%%%%

We first consider a generalized BSDE with the generator $\gen :\Omega \times [0,T]\times \RR \to \RR$ and the terminal
condition $\xi_T\in L^2(\cF_T)$. The mapping $\gen :\Omega\times[0,T]\times\RR\rightarrow \RR$ is
assumed to be $\mathcal{P}(\FF) \otimes\mathcal{B}(\RR)$-measurable where
$\mathcal{P}(\FF)$ is the $\sigma$-field of predictable sets in $\Omega\times[0,T]$.

\bd \label{ndef1.3}
A pair $(Y,Z)$ is a {\it solution to the generalized BSDE} with data $(M,\newA,\gen ,\xi_T)$ if
$Y$ is an $\FF$-adapted, c\`adl\`ag process, $Z$ is an $\FF$-predictable process, and the following equality holds,
for every $t \in [0,T]$,
\begin{align} \label{eq1.12m}
Y_t=\xi_T-\int_{]t,T]}Z_s\,dM_s+\int_{]t,T]} \gen (s,Y_s)\,d\newA_s
\end{align}
where the It\^o integral $\int_{]0,t]}Z_s\,dM_s$ is a martingale and the integral $\int_0^t \gen (s,Y_s)\,d\newA_s$
is an $\FF$-adapted, continuous process of finite variation.
\ed

% Notice that the process $Y$ in \eqref{eq1.12m} is a local supermartingale if the mapping $\gen $ is nonnegative.

We will frequently examine the existence and uniqueness of a solution $(Y,Z)$ to the BSDE \eqref{eq1.12m} in the product space $\cS^2(\FF) \times \cH^2(M)$. Then we say that a solution to the GBSDE \eqref{eq1.12m} is unique in $\cS^2(\FF) \times \cH^2(M)$ if for any two solutions $(Y',Z')$ and $(Y'',Z'')$ belonging to $\cS^2(\FF) \times \cH^2(M)$ we have that $||Y'-Y''||_{\cS^2(\FF)}=0$ and $||Z'-Z''||_{\cH^2(M)}=0$.

We start by establishing in Section \ref{nsec2.1} the comparison property for a GBSDE and we obtain in Section \ref{nsec2.2}
some {\it a priori} estimates when the driver $A$ is bounded. This allows us to study in Section \ref{nsec2.3} the well-posedness results for the GBSDE \eqref{eq1.12m} in the space $\cS^2(\FF) \times \cH^2(M)$ when the driver $A$ is either bounded or square-integrable. In Section \ref{nsec3.4}, for a given $\FF$-optional and bounded process $\eta$, we examine the sequence $(Y^n)_{n\in\NN}$ of solutions to the GBSDE
penalized GBSDEs
\begin{align*} % \label{veq6.9v}
Y^n_{\tau}= \zeta_T-\int_{\rrb\tau,T\rrb} Z^n_s\,dM_s+\int_{\rrb\tau,T\rrb} n(Y^n_s-\eta_s)^+\,d\newA_s
\end{align*}
and we show that the limit $Y:=\lim_{\,n\to\infty}Y^n$ can be interpreted as the value process of a particular optimal stopping problem.

%\bd \label{ndef1.3}
%A pair $(Y,Z)\in \cS^2(\FF) \times \cH^2(M)$ is a solution to the {\it generalized BSDE} with data $(M,\newA,\gen ,\xi_T)$ if
%the process $Y$ is c\`adl\`ag and the following equality holds,
%for every $t \in [0,T]$,
%\begin{align} \label{eq1.12m}
%Y_t=\xi_T-\int_{]t,T]}Z_s\,dM_s+\int_{]t,T]} \gen (s,Y_s)\,d\newA_s.
%\end{align}
%We say that a solution to the GBSDE \eqref{eq1.12m} is unique if for any two solutions $(Y',Z')$ and $(Y'',Z'')$ we have $||Y'-Y''||_{\cS^2(\FF)}=0$ and
% $||Z'-Z''||_{\cH^2(M)}=0$.
%\ed

%%%%%%%%%%%%%%%%%%%%%%%%%%%%%%%%%%%%%%%%%%%%%%%%%%%%%%%%%%%%%%%%%%%%%%%%%%%%%%%%%%%%%%%%%%%
\subsection{Comparison Theorem for a Generalized BSDE}           \label{nsec2.1}
%%%%%%%%%%%%%%%%%%%%%%%%%%%%%%%%%%%%%%%%%%%%%%%%%%%%%%%%%%%%%%%%%%%%%%%%%%%%%%%%%%%%%%%%%%%

Our first goal is to establish the comparison property of solutions to \eqref{eq1.12m} and thus we work under a temporary
assumption that a solution to the GBSDE \eqref{eq1.12m} exists but is not necessarily unique. For brevity, we denote
\begin{align*}
\Gen_t(Y):=\int_{]0,t]} \gen (s,Y_s)\,d\newA_s ,\quad H_t(\wtY):=\int_{]0,t]} h(s,\wtY_s)\,d\newA_s.
\end{align*}

\bp \label{npro1.1}
Let the mappings $\gen ,h:\Omega \times [0,T]\times \RR \to \RR$ be such that the GBSDEs
\begin{align*}
Y_t=\xi_T-\int_{]t,T]}Z_s\,dM_s+\int_{]t,T]} \gen (s,Y_s)\,d\newA_s=\xi_T-\int_{]t,T]}Z_s\,dM_s+\Gen_T(Y)-\Gen_t(Y)
\end{align*}
and
\begin{align*}
\wtY_t=\wtxi_T-\int_{]t,T]}\wtZ_s\,dM_s+\int_{]t,T]} h(s,\wtY_s)\,d\newA_s=\wtxi_T-\int_{]t,T]}\wtZ_s\,dM_s+H_T(\wtY)-H_t(\wtY)
\end{align*}
have solutions $(Y,Z)$ and $(\wtY,\wtZ)$. % in $\cS^2(\FF)\times \cH^2(M)$.
Assume that $\xi_T \ge \wtxi_T $ and the
functions $\gen$ and $h$ satisfy the following conditions: \\
(i) $\gen (\omega,t,y)\ge h(\omega,t,y)$ for every $(\omega,t,y)\in \Omega \times [0,T] \times \RR$, \\
(ii) $\gen (\omega,t,\cdot )$ % or $h(\omega,t,\cdot)$
 is a nonincreasing function for every $(\omega,t) \in \Omega \times [0,T]$. \\
Then the inequality $Y\ge \wtY $ is valid, that is, $\PP ( Y_t\geq \wtY_t,\,\forall\, t\in [0,T])=1$.
\ep

\begin{proof}
We adapt the proof of Lemma 8.3 in Peng \cite{P2004}. It is important to recall that the processes $Y$ and $\wtY$ are c\`adl\`ag.
For a fixed $\varepsilon>0$, we define the $\FF$-stopping time $\tau^{\varepsilon}:=\inf\,\{t\ge 0 : Y_t \leq \wtY_t - \varepsilon\}$ where, by convention, $\inf \emptyset = T$. Note that if for all $\varepsilon>0$ we have $\PP (\tau^{\varepsilon}=T)=1$, then $Y\ge \wtY$. Indeed, for $\varepsilon_n=\frac{1}{n}$, there exists $D_n$ such that $\PP(D_n)=0$ and $\tau^{\varepsilon_n}(\omega) =T$ for  $\omega\notin D_n$. Thus for  $\omega\notin D_n$, from the definition of $\tau^{\varepsilon}$ and the equality $\tau^{\varepsilon_n}(\omega)=T$,  we obtain $Y_t \geq \wtY_t- \varepsilon_n$ for every $t\in[0,T[$. Then for $D:=\cup_{n=1}^{\infty}D_n$, we have $\PP(D)=0$ and for $\omega\notin D$, we have  $Y_t \geq \wtY_t$ for every $t\in [0,T[$. Since $Y_T=\xi_T \geq \wtxi_T =\wtY_T$ we conclude that $Y\ge \wtY$.

We now argue by contradiction. If the inequality $Y\geq \wtY$ does not hold then, by the first step, there exists
$\varepsilon >0$ such that $\PP (E)>0$ where $E:=\{\tau^{\varepsilon} <T\} \in \cF_{\tau^{\varepsilon}}$.
We fix $\varepsilon$ and we define $\tau:=\tau^{\varepsilon}\I_E + T\I_{E^c}$ and $\nu:=\inf\,\{t\ge \tau : Y_t \geq  \wtY_t\}$ so that $\nu\leq T$ since, by assumption, $Y_T \geq \wtY_T$.  Since $Y$ and $\wtY$ are c\`adl\`ag processes it is clear that $Y_{\tau} < \wtY_{\tau}$ on $E$, the interval $\llb \tau,\nu \llb$ is nonempty on $E$, and the inequality $Y_{\nu} \geq \wtY_{\nu}$ is valid.

For brevity, let us write $U:=G(Y)$ and $\wtU:=H(\wtY)$ so that the $\FF$-adapted, continuous process $\ov{U} :=U-\wtU$ satisfies
\begin{align*}
\ov{U}_t& =\Gen_t(Y)-H_t(\wtY)=(\Gen_t(Y)-\Gen_t(\wtY))+(\Gen_t(\wtY)-H_t(\wtY)) \\
&= \int_{]0,t]} (g(s,Y_s)-g(s,\wtY_s))\, dA_s +  \int_{]0,t]} (g(s,\wtY_s)-h(s,\wtY_s))\, dA_s .
\end{align*}
We deduce that $\I_E \ov{U}$ is a continuous and nondecreasing process on $\llb\tau,\nu\rrb$ since  $\gen (\omega,t,y) \ge h(\omega,t,y)$ for all
$(\omega ,t,y) \in \Omega \times [0,T]\times \mathbb{R}$, the inequality $Y < \wtY$ holds on $\llb\tau,\nu\llb$, and for
every $(\omega ,t) \in \Omega \times [0,T]$ the function $\gen (\omega,t,\cdot )$ is nonincreasing.
We observe that on $E$
\begin{align*}
Y_t-\wtY_t =Y_{\nu}-\wtY_{\nu} -\int_{\rrb t,\nu \rrb}(Z_s-\wtZ_s) \,dM_s +\ov{U}_{\nu}-\ov{U}_t .
\end{align*}
Therefore, the process $Y-\wtY$ is a supermartingale (hence also strong supermartingale) on $\llb\tau,\nu\rrb$ and
$Y_{\nu}- \wtY_{\nu}\geq 0$. Consequently, $Y_{\tau}-\wtY_{\tau} %\geq \EE [ Y_{\nu}- \wtY_{\nu}\,|\,\cF_{\tau}]
\geq 0$ and thus $Y_{\tau}\geq \wtY_{\tau}$ on $E$, which leads to a contradiction since $Y_{\tau} < \wtY_{\tau}$ on $E \in \cF_{\tau}$ and $\PP (E)>0$.
\end{proof}

When using the penalization method to show the existence of a solution to the reflected GBSDE with a lower obstacle $\xi$, one
may employ the mapping $f:\Omega \times [0,T]\times \RR \to \RR_+$, which is given by $f(t,y) =(\xi_t-y)^+$ where $\xi$ is a
predetermined $\FF$-optional process. The following corollary to Proposition \ref{npro1.1} will be useful in the proof
of the penalization result (see Section \ref{nsec3.4}).

\bcor  \label{ncor1.1}
Let the mapping $f:\Omega \times [0,T]\times \RR \to \RR_+$ be nonnegative and such that $f(t,\cdot)$ is nonincreasing. For every $n\in\NN$, let the pair
$(Y^n,Z^n)$ be a %unique
solution to the GBSDE
\begin{align} \label{eq1.12v}
Y^n_t=\xi_T-\int_{]t,T]}Z^n_s\,dM_s+\int_{]t,T]} n f(s,Y^n_s)\,d\newA_s.
\end{align}
Then $Y^{n+1}\geq Y^n$ for every $n\in\NN$, that is, $\PP (Y^{n+1}_t\geq Y^n_t,\, \forall \, t \in [0,T])=1$.
\ecor

\begin{proof}
It suffices to apply Proposition \ref{npro1.1} to the GBSDEs with generators $\gen (t,y)=(n+1)f(t,y)$ and $h(t,y)=n f(t,y)$.
\end{proof}

%%%%%%%%%%%%%%%%%%%%%%%%%%%%%%%%%%%%%%%%%%%%%%%%%%%%%%%%%%%%%%%%%%%%%%%%%%%%%%%%%%%%%%%%%%%%%%%%
\subsection{A Priori Estimates for Solutions to a Generalized BSDE}     \label{nsec2.2}
%%%%%%%%%%%%%%%%%%%%%%%%%%%%%%%%%%%%%%%%%%%%%%%%%%%%%%%%%%%%%%%%%%%%%%%%%%%%%%%%%%%%%%%%%%%%%%%%

In addition to Assumption \ref{nass1.1}, we make the following standing assumption on the generator $\gen$ of the GBSDE \eqref{eq1.12m}.

\bhyp \label{nass1.2}
Let the mapping $\gen :\Omega \times [0,T]\times \RR \to \RR$ satisfy, for every $t\in[0,T]$,
\begin{align*}
&(i)  \ \ \EE\bigg[\int_{]0,T]} |\gen (t,0)|^2\, d\newA_t\bigg]< +\infty ,
&(ii) \ \ |\gen (t,y)-\gen (t,y')| \leq L|y-y'|,\ \forall\, y, y' \in \RR.
\end{align*}
\ehyp

Our first goal is to establish a useful {\it a priori} estimate for a postulated solution in $\cS^2(\FF)\times \cH^2(M)$ to the GBSDE \eqref{eq1.12m}.

\bp \label{npro1.2}
Let Assumptions \ref{nass1.1} and \ref{nass1.2} hold with a bounded process $A$. If $(Y,Z)\in \cS^2(\FF)\times \cH^2(M)$ is a solution to the GBSDE \eqref{eq1.12m}, then
there exists a constant $c>0$ such that for every $\alpha, \beta>0$ such that $\beta >2L+\alpha^{-1}$ we have
\begin{align} \label{neq1.14}
\EE\bigg[\sup_{t\in[0,T]}\big[ \ee_t|Y_t|^2\big]+\int_{]0,T]} \ee_s|Z_s|^2\,d[M]_s\bigg] \leq c \,\EE \bigg[ e_T|\yxi_T|^2+ \alpha\,\int_{]0,T]} \ee_s|g(s,0)|^2\,d\newA_s \bigg]
\end{align}
where $\ee_t := e^{\beta \newA_t}$.
\ep

\begin{proof}
Note that $d\ee_t := \beta \ee_t \, d\newA_t$ and $1 \leq \ee_t \leq e^{\beta c_{\newA}}$ for every $t\in[0,T]$ since $\beta >0$ and thus $e$ is a bounded process. We first establish the {\it a priori} estimates for a postulated solution $(Y,Z)\in\cS^2(\FF)\times\cH^2(M)$ to the BSDE \eqref{eq1.12m}.
By applying the It\^o formula to $\ee_t|Y_t|^2$ and the Young inequality with $\alpha>0$, we obtain
\begin{align*}
&\ee_t|Y_t|^2+\int_{]t,T]}\ee_s|Z_s|^2\,d[M]_s \\
&\leq e_T|\yxi_T|^2-\beta\int_{]t,T]}\ee_s|Y_s|^2\,d\newA_s+2\int_{]t,T]}\ee_sY_s\big(g(s,0)+L|Y_s|\big)\,dA_s-2\int_{]t,T]}\ee_sY_{s-}Z_s\,dM_s \\
&\leq e_T|\yxi_T|^2+\alpha\,\int_{]t,T]}\ee_s|g(s,0)|^2\,d\newA_s+\big(2L+\alpha^{-1}-\beta\big)\,\int_{]t,T]}\ee_s|Y_s|^2\,d\newA_s
-2\int_{]t,T]}\ee_sY_{s-}Z_s\,dM_s 
\end{align*}
and thus
\begin{align}\label{eq1.16a}
\ee_t|Y_t|^2+\int_{]t,T]}\ee_s|Z_s|^2\,d[M]_s \leq e_T|\yxi_T|^2+\alpha\,\int_{]t,T]}\ee_s|g(s,0)|^2\,d\newA_s
-2\int_{]t,T]}\ee_sY_{s-}Z_s\,dM_s
\end{align}
for any $\alpha,\beta>0$ such that $\beta > 2L+\alpha^{-1}$. By taking the supremum and expectation in \eqref{eq1.16a} we obtain
\begin{align} \label{eq1.16c}
\EE\bigg[\sup_{t\in [0,T]} \ee_t|Y_t|^2\bigg]\leq \EE \bigg[e_T|\yxi_T|^2+ \alpha\,\int_{]0,T]} \ee_s|g(s,0)|^2\, d\newA_s\bigg]
 +2\,\EE\bigg[\esssup_{t\in [0,T]}\Big|\int_{]0,t]}\ee_sY_{s-}Z_s\,dM_s\Big|\bigg].
\end{align}
Since $A$ is a bounded process an application of the Burkholder--Davis--Gundy inequality with $p=1$ gives % (that is, the Davis inequality) gives
% (see, e.g., Theorem 10.24 in He et al. \cite{HWY1992} or Theorem 11.5.5 in Cohen and Elliott \cite{CE2015} with $p=1$)
\begin{align} \label{eq1.16d}
&2\,\EE\bigg[\sup_{t\in [0,T]} \Big|\int_{]0,t]}\ee_sY_{s-}Z_s\,dM_s\Big|\bigg]
\leq 2c_1\,\EE \bigg[\bigg(\int_{]0,T]} e^2_s|Y_{s-}|^2|Z_s|^2\,d[M]_s\bigg)^{1/2}\bigg] \nonumber \\
&\leq \EE \bigg[\bigg(\frac{1}{2}\sup_{t\in [0,T]}\ee_t|Y_t|^2\bigg)^{1/2}\bigg(8c_1^2\int_{]0,T]}\ee_s|Z_s|^2\,d[M]_s\bigg)^{1/2}\bigg]\\
&\leq \frac{1}{4}\,\EE \bigg[\sup_{t\in [0,T]}\ee_t|Y_t|^2\bigg]
+4c_1^2\,\EE\bigg[\int_{]0,T]} \ee_s|Z_s|^2\,d[M]_s\bigg] \nonumber
\end{align}
where the constant $c_1$ % $c_1=2\sqrt{6}$
is independent of $\alpha,\beta,L$ and the last inequality holds since $2ab\leq a^2+b^2$ for all real numbers $a,b$.
Furthermore, by taking the expectation in \eqref{eq1.16a} and using the martingale property of the integral with respect to $M$, we obtain
\begin{align} \label{eq1.16b}
\EE\bigg[ \int_{]0,T]}\ee_s|Z_s|^2\,d[M]_s \bigg]\leq	\EE \bigg[ e_T|\yxi_T|^2+ \alpha\,\int_{]0,T]} \ee_s|g(s,0)|^2\,d\newA_s \bigg]<\infty
\end{align}
where the second inequality follows from Assumption \ref{nass1.2} (i) and the boundedness of $A$.
By combining \eqref{eq1.16c}, \eqref{eq1.16d} and \eqref{eq1.16b}, we obtain
\begin{align*}
\frac{3}{4}\,\EE\bigg[\sup_{t\in[0,T]} \ee_t|Y_t|^2 \bigg]\leq (1+4c_1^2)\, \EE \bigg[ e_T|\yxi_T|^2+ \alpha\,\int_{]0,T]} \ee_s|g(s,0)|^2\,d\newA_s \bigg]
\end{align*}
and thus
\begin{align*}
\EE\bigg[\sup_{t\in[0,T]} \ee_t|Y_t|^2 +\int_{]0,T]} \ee_s|Z_s|^2\,d[M]_s\bigg] \leq c\, \EE \bigg[ e_T|\yxi_T|^2+ \alpha\,\int_{]0,T]} \ee_s|g(s,0)|^2\,d\newA_s \bigg]
\end{align*}
where the constant $c:=\frac{7+16c_1^2}{3}>0$ is independent of $\alpha,\beta$ and $L$. 
% We stress that the boundedness of $\newA$ was not used in the proof.
\end{proof}

The next results deals with the stability of solutions to a GBSDE with respect to the terminal condition $\xi_T$ and generator $g$.
Let us denote $\wh Y=Y^1-Y^2,\, \wh Z=Z^1-Z^2$ and $\wh g_t = g^1(t,Y^2_t)- g^2(t,Y^2_t)$. As in Proposition \ref{npro1.2}, we denote 
$\ee_t:=e^{\beta \newA_t}$.

\bp  \label{npro1.2x}
Let Assumptions \ref{nass1.1} and \ref{nass1.2} hold for $g^i,\, i=1,2$ and the nondecreasing, continuous process $\newA$ be bounded on $[0,T]$ so that $\newA_T\leq c_{\newA}$ for some constant $c_{\newA}$. If $(Y^i,Z^i)\in \cS^2(\FF)\times \cH^2(M)$ is a solution to the GBSDE with data  $(M,\newA,\gen^i ,\xi^i_T)$ for $i=1,2$, then for every $\alpha, \beta>0$ such that $\beta >2L+\alpha^{-1}$, we have, for every $t\in [0,T]$,
\begin{align}  \label{neq1.14x}
\ee_t|\wh Y_t|^2 \leq \,\EE \bigg[\ee_T|\wh\xi_T|^2+\alpha\int_{]t,T]} \ee_s|\wh\gen_s|^2\,d\newA_s\,\Big|\,\cF_t\bigg].
\end{align}
\ep

\begin{proof}
An application of the It\^o formula to $\ee_t|\whY_t|^2$ from $[t,T]$ and the Young inequality with $\alpha>0$ gives
\begin{align*} %\label{neq1.15x}
&\ee_t|\whY_t|^2+\int_{]t,T]}\ee_s|\whZ_s|^2\,d[M]_s \nonumber\\
&\leq e_T|\wh\yxi_T|^2-\beta\int_{]t,T]}\ee_s|\whY_s|^2\,d\newA_s+2\int_{]t,T]}\ee_s\wh Y_s\big(L|\whY_s|+|\wh g_s|\big)\,dA_s-2\int_{]t,T]}\ee_s\wh Y_{s-}\wh Z_s\,dM_s \\
&\leq e_T|\wh\yxi_T|^2+\alpha\,\int_{]t,T]}\ee_s|\wh g_s|^2\,d\newA_s+\big(2L+\alpha^{-1}-\beta\big)\,\int_{]t,T]}\ee_s|\wh Y_s|^2\,d\newA_s
-2\int_{]t,T]}\ee_s\whY_{s-}\wh Z_s\,dM_s \nonumber
\end{align*}
and thus for any $\alpha,\beta>0$ such that $\beta>2L+\alpha^{-1}$
\begin{align*}  %\label{neq1.15x}
\ee_t|\wh Y_t|^2\leq  e_T|\wh\yxi_T|^2+\alpha\,\int_{]t,T]}\ee_s|\wh g_s|^2\,d\newA_s-2\int_{]t,T]}\ee_s\whY_{s-}\whZ_s\,dM_s.
\end{align*}
 To complete the proof it suffices to take the conditional expectation with respect to $\cF_t$.
\end{proof}

%%%%%%%%%%%%%%%%%%%%%%%%%%%%%%%%%%%%%%%%%%%%%%%%%%%%%%%%%%%%%%%%%%%%%%%%%%%%%%%%%%%%%%%%%%%%%%%%
\subsection{Existence of a Solution to a Generalized BSDE}      \label{nsec2.3}
%%%%%%%%%%%%%%%%%%%%%%%%%%%%%%%%%%%%%%%%%%%%%%%%%%%%%%%%%%%%%%%%%%%%%%%%%%%%%%%%%%%%%%%%%%%%%%%%

The next two propositions deal with the existence and uniqueness of a solution to the GBSDE \eqref{eq1.12m} with a Lipschitz
continuous driver. We first present the existence result under an additional postulate that the driver $A$ is bounded.
It will be subsequently extended in Proposition \ref{npro1.4} to the case of Assumption \ref{nass1.1} (iii) where the driver is assumed to be square-integrable.

\bp \label{npro1.3}
Let Assumptions \ref{nass1.1} and \ref{nass1.2} hold with a bounded process $A$. Then the GBSDE
\begin{align} \label{eq1.15a}
Y_t=\yxi_T-\int_{]t,T]}Z_s\,dM_s+\int_{]t,T]} \gen (s,Y_s)\,d\newA_s
\end{align}
has a unique solution $(Y,Z)$ in the space $\cS^2(\FF)\times \cH^2(M)$.
\ep

\begin{proof}
To prove the existence and uniqueness of a solution to the GBSDE \eqref{eq1.15a}, we use the standard method based on the Banach fixed point theorem. For any $Y\in \cS^2(\FF)$ and $\beta > 0$, we define $||Y||_{\cS^2_\beta(\FF)}^2:=\EE\big[\sup_{t \in [0,T]}\ee_t |Y_t|^2\big]$ where $\ee_t= e^{\beta \newA_t}$ and we observe that the norms $||\cdot||_{\cS^2(\FF)}$ and $||\cdot||_{\cS^2_\beta(\FF)}$ are equivalent on $\cS^2(\FF)$. We denote by $\cS_{\beta}^2(\FF)$ the space $\cS^2(\FF)$ endowed with the norm $||\cdot||_{\cS^2_\beta(\FF)}$.
Let the mapping $\Phi : \cS_{\beta}^2(\FF) \to \cS_{\beta}^2(\FF)$ be defined as follows: for any given $\ww \in \cS_{\beta}^2(\FF)$ we set $\Phi(\ww):=Y^{\ww}$ where the pair $(Y^{\ww},Z^{\ww})$ is a unique solution to the GBSDE
\begin{align} \label{neq1.17}
Y^{\ww}_t=\yxi_T-\int_{]t,T]}Z^{\ww}_s\,dM_s+\int_{]t,T]} g(s,\ww_s)\,d\newA_s
\end{align}
and the fixed generator is clearly independent of $Z^{\ww}$. To show that a unique solution to \eqref{neq1.17} exists, we observe that the process $\Theta$, which given by
\begin{align} \label{eq1.17a}
\Theta_t:=Y^{\ww}_t+\int_{]0,t]} g(s,\ww_s)\,d\newA_s=Y^{\ww}_0+\int_{]0,t]}Z^{\ww}_s\,dM_s,
\end{align}
is a martingale where the existence of a process $Z^{\ww}\in\cH^2(M)$ follows from the postulated predictable representation property of $M$.
Hence the unique solution $(Y^{\ww},Z^{\ww})\in \cS^2_{\beta}(\FF)\times \cH^2(M)$ to \eqref{neq1.17} is determined by
\begin{align*}
Y^{\ww}_t=\EE\bigg[\yxi_T + \int_{]t,T]} g(s,\ww_s)\,d\newA_s\,\Big|\,\cF_t \bigg]=:\Phi (\ww).
\end{align*}

Our next goal is to demonstrate that there exists a unique process $\wh \ww \in \cS_{\beta}^2(\FF) $ such that $\Phi (\wh \ww)=\wh \ww$.
Then the corresponding process $\wh z \in \cH^2(M)$ can be found from equality \eqref{eq1.17a} combined with the predictable representation property of $M$, that is, from the equality
\begin{align*}
\wh \ww_t +\int_{]0,t]} g(s,\wh \ww_s)\,d\newA_s=\wh \ww_0+\int_{]0,t]}\wh z_s\,dM_s.
\end{align*}
It is clear that it suffices to show that the mapping $\Phi:\cS_{\beta}^2(\FF)\to\cS_{\beta}^2(\FF)$ is a contraction when $\beta$ is sufficiently large.
To this end, we take $\ww', \ww''\in \cS_{\beta}^2(\FF)$ and denote $Y^{\ww'}=\Phi(\ww')$ and $Y^{\ww''}=\Phi(\ww'')$. For the simplicity of notation, we write $\yyy:=Y^{\ww'}-Y^{\ww''}=\Phi(\ww')-\Phi(\ww''),\, \zzz :=Z^{\ww'}-Z^{\ww''}$ and $\ww:= \ww'- \ww''$.
It is clear from \eqref{neq1.17} that $\yyy$ satisfies the GBSDE
\begin{align*}
\yyy_t=-\int_{]t,T]}\zzz_s\,dM_s+\int_{]t,T]} \big(g(s,\ww'_s)-g(s,\ww''_s)\big)\,d\newA_s
\end{align*}
where $|g(s,\ww'_s)-g(s,\ww''_s)|\leq L|\ww_s|$ since $g(s,\cdot)$ is a Lipschitz continuous function with constant $L$.

The foregoing computations are similar to those used in the proof of Proposition \ref{npro1.2}, so we only sketch some significant steps.
By applying the It\^o formula to $\ee_t|\yyy_t|^2$ and subsequently the Young inequality with $\alpha>0$, we obtain
\begin{align} \label{neq1.18}
&\ee_t|\yyy_t|^2+\int_{]t,T]}\ee_s|\zzz_s|^2\,d[M]_s \leq -\beta\int_{]t,T]}\ee_s|\yyy_s|^2\,d\newA_s +2L \int_{]t,T]} \ee_s \yyy_s\,|\ww_s|d\newA_s -2\int_{]t,T]}\ee_s\yyy_{s-}\zzz_s\,dM_s \nonumber\\
&\leq -\beta\,\int_{]t,T]}\ee_s|\yyy_s|^2\,d\newA_s + \frac{1}{\alpha}\,\int_{]t,T]}\ee_s|\yyy_s|^2\,d\newA_s
+ \alpha L^2\,\int_{]t,T]} \ee_s|\ww_s|^2\, d\newA_s-2\int_{]t,T]}\ee_s\yyy_{s-}\zzz_s\,dM_s\\
&\leq \alpha L^2\,\int_{]t,T]} \ee_s|\ww_s|^2\, d\newA_s-2\int_{]t,T]}\ee_s\yyy_{s-}\zzz_s\,dM_s \nonumber
\end{align}
where we have assumed that $\beta > \alpha^{-1} $. For brevity, we denote
\begin{align*}
J(w):=\alpha L^2\,\int_{]0,T]} \ee_s|\ww_s|^2\,d\newA_s.
\end{align*}
By setting $t=0$ in \eqref{neq1.18} and taking the expectation, we obtain
\begin{align}\label{neq1.19}
\EE\bigg[\int_{]0,T]}\ee_s|\zzz_s|^2\,d[M]_s\bigg]\leq \EE[J(w)].
\end{align}
Furthermore, it follows from \eqref{neq1.17} that, for every $t\in [0,T]$,
\begin{align*}
\ee_t|\yyy_t|^2&\leq J(w) -2\int_{]t,T]}\ee_s\yyy_{s-}\zzz_s\,dM_s \leq J(w) +2\Big|\int_{]t,T]}\ee_s\yyy_{s-}\zzz_s\,dM_s\Big|
\end{align*}
and thus, by taking the supremum over $t\in [0,T]$ and the expectation on both sides, we get
\begin{align} \label{neq1.20}
\EE\bigg[\sup_{t\in [0,T]} \ee_t|\yyy_t|^2\bigg]&\leq \EE [J(w)] +2\,\EE\bigg[\esssup_{t\in [0,T]}\Big|\int_{]t,T]}\ee_s\yyy_{s-}\zzz_s\,dM_s\Big|\bigg]\\
&\leq \EE [J(w)]+ \frac{1}{4}\,\EE \bigg[\sup_{t\in [0,T]}\ee_t|\yyy_t|^2\bigg]
+4c_1^2\,\EE\bigg[\int_{]0,T]} \ee_s|\zzz_s|^2\,d[M]_s\bigg] \nonumber
\end{align}
where we have used the Burkholder--Davis--Gundy inequality similarly as in \eqref{eq1.16d} and thus the constant $c_1$ % $c_1=2\sqrt{6}$
is independent of $L, \alpha$ and $\beta$. By combining \eqref{neq1.19} and \eqref{neq1.20}, we obtain
\begin{align*}
\frac{3}{4}\,\EE\bigg[\sup_{t\in[0,T]} \ee_t|\yyy_t|^2 \bigg]\leq (1+4c_1^2)\,\EE [J(w)],
\end{align*}
which in turn implies that
\begin{align*}
&\frac{3}{4}\, \EE\bigg[\sup_{t\in [0,T]} \ee_t|\yyy_t|^2 \bigg]
\leq \alpha L^2(1+4c_1^2)\,\EE\bigg[\int_{]0,T]} \ee_s|\ww_s|^2 d\newA_s\bigg] \\
& \leq \alpha L^2(1+4c_1^2)\,\EE\bigg[\newA_T \sup_{t\in [0,T]} \ee_t|\ww_t|^2\bigg]
\leq \alpha  c_{\newA} L^2 (1+4c_1^2)\, \EE\bigg[ \sup_{t\in [0,T]} \ee_t|\ww_t|^2\bigg]
%& \leq n^2\alpha(1+4c_1^2)\bigg[\EE\bigg(\sup_{t\in [0,T]} \ee_t|\ww_t|^2\bigg)^2 \bigg]^{1/2}\bigg[\EE\big[(\newA_T-\newA_0)^2\big]\bigg]^{1/2}
\end{align*}
since $\newA_T \leq c_{\newA}$. Consequently, for all $w',w''\in\cS_{\beta}^2(\FF)$,
\begin{align*}
||\Phi(\ww')-\Phi(\ww'')||_{\cS^2_\beta(\FF)}^2\leq\frac{4}{3}\alpha c_{\newA} L^2 (1+4c_1^2)\,||w'-w''||_{\cS^2_\beta(\FF)}^2 = \gamma ||w'-w''||_{\cS^2_\beta(\FF)}^2.
\end{align*}
We conclude that $\Phi $ is a contraction when $\beta > \alpha^{-1}$ and $\alpha>0$ is such that $\gamma <1$, that is, when $\alpha < \frac{3}{4}c^{-1}_{\newA}{L}^{-2}(1+4c_1^2)^{-1}$ and $\beta > \alpha^{-1}$. Then, from the Banach theorem, there exists a unique fixed point of the mapping $\Phi: \cS_{\beta}^2(\FF) \to \cS_{\beta}^2(\FF)$. This shows that the GBSDE \eqref{eq1.15a} has a unique solution when the inequality $\newA_T\leq c_{\newA}$ holds.
\end{proof}

In the next proposition we relax the assumption that the driver $A$ is bounded and we work instead under Assumption \ref{nass1.1} (iii) that $A$ 
is a square-integrable process. However, since the method of proof hinges on the comparison property we need to assume that the generator is nonnegative.

\bp \label{npro1.4}
Let Assumptions \ref{nass1.1} and \ref{nass1.2} hold and the function $g(\omega,t,\cdot)$ be nonnegative and nonincreasing for every $(\omega,t)\in \Omega \times [0,T]$. Then the GBSDE
\begin{align} \label{eq1.15b}
Y_t=\yxi_T-\int_{]t,T]}Z_s\,dM_s+\int_{]t,T]} g(s,Y_s)\,d\newA_s
\end{align}
has a unique solution $(Y,Z)$ in $\cS^2(\FF)\times \cH^2(M)$.
\ep

\begin{proof} Define the sequence $\newA^n := \newA \wedge n$ of continuous, nondecreasing, bounded processes so that $c_{\newA^n}=n$ for every $n\in\NN$. We will use the comparison property for solutions to the GBSDE in order to obtain a solution
to \eqref{eq1.15b} with $\newA$ as a limit of solutions to  \eqref{eq1.15a} with $\newA$ replaced by $\newA^n$.

We already know from Proposition \ref{npro1.3} that for every $n\in\NN$ the GBSDE
\begin{align} \label{eq1.15c}
Y^{n}_t=\yxi_T-\int_{]t,T]}Z^{n}_s\,dM_s+\int_{]t,T]}g(s,Y_s^n)\,d\newA^n_s
\end{align}
has a unique solution $(Y^{n},Z^{n})\in \cS^2(\FF)\times \cH^2(M)$.
To complete the proof, we will show the following assertions: \\
(a) the sequence $(Y^{n})_{n\in\NN}$ is increasing, in the sense that $Y^{n+1}\geq Y^{n}$ for every $n\in\NN$; \\
(b) there exists a constant $c$ independent of $n$ such that
\begin{align} \label{eq1.20c}
\EE\bigg[\sup_{t\in[0,T]}\ee_t|Y_t^n|^2+\int_{]0,T]} \ee_s|Z_s^n|^2\,d[M]_s\bigg]\leq c\,\EE \bigg[e_T|\yxi_T|^2
+\alpha\,\int_{]0,T]} \ee_s|g(s,0)|^2\,d\newA^n_s \bigg],
\end{align}
(c) the limit $(Y,Z)=\lim_{\,n\to \infty} (Y^n,Z^n)$ in $\cS^2(\FF)\times \cH^2(M)$ is a unique solution to the GBSDE \eqref{eq1.15a}.

\vskip 5 pt
\noindent (a) To establish the claimed comparison property $Y^{n+1}\geq Y^n$, we fix $n$ and we apply the method from the proof of
Proposition \ref{npro1.1} to the pair $Y$ and $\wtY$ where $Y:=Y^{n+1}$ and $\wtY := Y^n$. We set $\xi = \eta = \yxi_T$ and $\gen =h$ and
we define
\begin{align*}
U_t :=\Gen_t(Y)=\int_{]0,t]} \gen (s,Y_s)\,d\newA^{n+1}_s, \quad \wtU_t:= H_t(\wtY)=\int_{]0,t]} \gen (s,\wtY_s)\,d\newA^{n}_s.
\end{align*}
Then we have that
\begin{align*}
V_t&:= U_t-\wtU_t=\Gen_t(Y)-H_t(\wtY)=(\Gen_t(Y)-\Gen_t(\wtY))+(\Gen_t(\wtY)-H_t(\wtY)) \\
&= \int_{]0,t]} (g(s,Y_s)-g(s,\wtY_s))\, dA^{n+1}_s +  \int_{]0,t]} g(s,\wtY_s)\,d(\newA^{n+1}_s-\newA^{n}_s).
\end{align*}
Since the generator $g(\omega,t,y)$ is assumed to be nonnegative and nonincreasing in $y$ and the process $\newA^{n+1}-\newA^{n}$ is nondecreasing,
the continuous and $\FF$-adapted process $\I_E V$ is nondecreasing on $\llb\tau,\nu\rrb$ where the event $E$ and the stopping times $\tau$ and $\nu$ are defined as in the proof of Proposition \ref{npro1.1}. To complete the proof of the proposition, it now suffices to check that all other arguments from the proof of Proposition \ref{npro1.1} are still valid.

\vskip 5 pt
\noindent (b)  An inspection of the proof of Proposition \ref{npro1.2} shows that the constants $c=\frac{7+16c_1^2}{3}$ and $c_1$
are independent of the process $\newA$ and thus assertion (b) is valid.

\vskip 5 pt
\noindent (c) The uniqueness of a solution to \eqref{eq1.15b} has already been established in Proposition \ref{npro1.1} and thus it remains to show that  $(Y,Z)=\lim_{\,n\to \infty} (Y^n,Z^n)$ is well defined in $\cS^2(\FF)\times\cH^2(M)$ and is a solution to the GBSDE \eqref{eq1.15b}.

First, we note that the sequence $(Y^n)_{n\in\NN}$ converges increasingly to a process $Y$. From inequality \eqref{eq1.20c}, there exists a constant $c$ independent of $n$ such that $||Y^n||_{\cS_{\beta}^2(\FF)} \leq c$ since, by assumption, $\xi_T$ is square-integrable and $\EE\big[\int_{]0,T]} |g(s,0)|^2\,d\newA_s\big]<\infty $. Hence we deduce that $\lim_{\,n\to \infty}||Y^n-Y||_{\cS_{\beta}^2(\FF)}^2=0$, which means that $(Y^n)_{n\in\NN}$ converges to $Y$ in the space $\cS_{\beta}^2(\FF)$. Consequently, the sequence $(Y^n)_{n\in\NN}$ is a Cauchy sequence in $\cS_{\beta}^2(\FF)$.

Second, we will demonstrate the convergence of the sequence $(Z^n)_{n\in\NN}$. For every $n,m\in \NN$ such that $n\geq m$, we observe that
\begin{align*}
\kkk^{n,m}_t&:= \int_{]0,t]} g(s,Y^n_s)\,d\newA^n_s
-\int_{]0,t]} g(s,Y^m_s)\,d\newA^m_s  \\
&\leq  \int_{]0,t]} \big| g(s,Y^n_s)- g(s,Y^m_s)\big| \,d\newA^n_s \\
&\leq \int_{]0,t]} L|Y^n_s-Y^m_s|\,d\newA^n_s
% \leq \int_{]0,t]} c_{\varphi}|Y^{\varphi,n}_s-Y^{\varphi,m}_s|\,d\newA_s
= L \int_{]0,t]} |\yyy_s|\,d\newA^n_s
\end{align*}
where we denote $\yyy_t =\yyy^{n,m}_t:= Y^n_t-Y^m_t \geq 0$. Let us set $\zzz_t = \zzz^{n,m}_t := Z^n_t-Z^m_t$. By applying the It\^o formula to $\ee^n_t|\yyy_t|^2$ where $\ee^n_t:= e^{\beta \newA^n_t}$ and the Young inequality with $\alpha>0$, we obtain
\begin{align*} % \label{eq1.21b}
&\ee^n_t|\yyy_t|^2+\int_{]t,T]}\ee^n_s|\zzz_s|^2\,d[M]_s=-\beta\int_{]t,T]}\ee^n_s|\yyy_s|^2\,d\newA^n_s +2\int_{]t,T]} \ee^n_s \yyy_s\,d\kkk^{n,m}_s-2\int_{]t,T]}\ee^n_s\yyy_{s-}\zzz_s\,dM_s \nonumber\\
&\leq -\beta\,\int_{]t,T]}\ee^n_s|\yyy_s|^2\,d\newA^n_s + \frac{1}{\alpha}\,\int_{]t,T]}\ee^n_s|\yyy_s|^2\,d\newA^n_s
+ L^2\alpha\,\int_{]t,T]} \ee^n_s|\yyy_s|^2\, d\newA^n_s-2\int_{]t,T]}\ee^n_s\yyy_{s-}\zzz_s\,dM_s\\
&\leq \alpha L^2\,\int_{]t,T]} \ee^n_s|\yyy_s|^2\, d\newA^n_s-2\int_{]t,T]}\ee^n_s\yyy_{s-}\zzz_s\,dM_s \nonumber
\end{align*}
where we have assumed that $\beta > \alpha^{-1}$. By setting $t=0$ and taking the expectation, we obtain
\begin{align}\label{eq1.22b}
\EE\bigg[\int_{]0,T]}\ee^n_s|\zzz_s|^2\,d[M]_s\bigg]\leq L^2\,\EE\bigg[\int_{]0,T]} \ee^n_s|\yyy_s|^2\,d\newA^n_s\bigg].
\end{align}
Consequently,
\begin{align*}
&||Z^n-Z^m||_{\cH^2(M)}^2=\EE\bigg[\int_{]0,T]}|Z^n_s-Z^m_s|^2\,d[M]_s\bigg] \leq
\EE\bigg[\int_{]0,T]}\ee^n_s|Z^n_s-Z^m_s|^2\,d[M]_s\bigg] \\
& %\leq \alpha L^2\,\EE\bigg[\int_{]0,T]} \ee^n_s|Y^n-Y^m|^2\,d\newA^n_s\bigg]
\leq \alpha L^2\,\EE\bigg[\int_{]0,T]} \ee^n_s|Y^n-Y^m|^2\,d\newA^n_s\bigg]
=\alpha L^2 n e^{n\beta T} \,||Y^n-Y^m||_{\cS^2(\FF)}^2.
\end{align*}
Since $(Y^{n_k})_{k\in\NN}$ is a convergent sequence, there exists a further subsequence $(Z^{m_k})_{k\in\NN}$ such that
\begin{align*}
\sum_{k=1}^\infty ||Z^{m_{k+1}}-Z^{m_k}||_{\cH^2(M)} < \infty ,
\end{align*}
which shows that $(Z^{m_k})_{k\in\NN}$ is a Cauchy sequence in the space $\cH^2(M)$ so that it converges in $\cH^2(M)$ to a limit $Z$.
It is now easy to conclude that the pair $(Y,Z)\in \cS^2(\FF)\times\cH^2(M)$ is a unique solution to the GBSDE \eqref{eq1.15b}.
\end{proof}

It is worth noting that if the generator satisfies the assumption of Proposition \ref{npro1.4} and is bounded so that $|g|\leq c_g$ for
some constant $c_g$, then the sequence $(K^n)_{n\in\NN}$ of increasing processes given by
\begin{align*}
K^n_t:=\int_{]0,t]}\gen (s,Y^n_s)\,d\newA^n_s
\end{align*}
converges to the increasing process $K$, which equals
\begin{align*}
K_t:=\int_{]0,t]}\gen (s,Y_s)\,d\newA_s
\end{align*}
and the convergence is uniform in $t$, for almost all $\omega$. We have
\begin{align*}
&\big|K_t- K^n_t\big|=\bigg|\int_{]0,t]}\gen (s,Y_s)\,d\newA_s-\int_{]0,t]}\gen (s,Y^n_s)\,d\newA^n_s \bigg| \\
&\leq \bigg|\int_{]0,t]}\gen (s,Y_s)\,d\newA_s-\int_{]0,t]}\gen (s,Y^n_s)\,d\newA_s \bigg|
+\bigg|\int_{]0,t]}\gen (s,Y^n_s)\,d\newA_s-\int_{]0,t]}\gen (s,Y^n_s)\,d\newA^n_s \bigg|\\
&\leq L \int_{]0,t]}\big|Y_s -Y^n_s\big|\,d\newA_s+c_g L \big(\newA_T-\newA^n_T\big)
\end{align*}
and thus
\begin{align*}
\sup_{t\in [0,T]}\big|K_t-K^n_t\big|
&\leq L\int_{]0,T]} \big|Y_s-Y^n_s\big|\,d\newA_s+c_g L\big(\newA_T-\newA^n_T\big)\\
&\leq L \newA_T\sup_{t\in [0,T]}\big|Y_t-Y^n_t\big|+ c_gL\big(\newA_T-\newA^n_T\big),
\end{align*}
which entails that $\sup_{t\in [0,T]}|K_t-K^n_t|$ converges to 0 almost surely when $n$ tends to $\infty$.

\bcor \label{ncor1.2}
Let Assumptions \ref{nass1.1} and \ref{nass1.2} hold and $\eta$ be an $\FF$-optional and bounded process. Then for every $n \in \NN$ the GBSDE
\begin{align} \label{eq1.16}
Y^n_t=\yxi_T-\int_{]t,T]}Z^n_s\,dM_s+\int_{]t,T]} n(\eta_s-Y^n_s)^+\,d\newA_s
\end{align}
has a unique solution $(Y^n,Z^n)$ in $\cS^2(\FF)\times \cH^2(M)$ and the inequality $Y^{n+1} \geq Y^n$ is satisfied for every $n \in \NN$.
\ecor

\begin{proof}
The first assertion is an immediate consequence of Proposition \ref{npro1.4}. To show that $Y^{n+1}\geq Y^n$, we observe that it suffices to take $\gen (t,y):= (n+1)(\eta_t-y)^+$ and $h(t,y):= n(\eta_t-y)^+$ in Proposition \ref{npro1.1} (see also Corollary \ref{ncor1.1}).
\end{proof}

It is worth noting that since $Y^n  \leq Y^{n+1}$ there exists a process $Y$ such that $Y^n \nearrow Y$. Furthermore, by the monotone convergence theorem
and Proposition \ref{npro1.2}, if the process $A$ is bounded then there exists a constant $c$ such that
\begin{align*}
\EE \bigg[\sup_{0 \leq t \leq T} |Y_t|^2\bigg] \leq c.
\end{align*}

%%%%%%%%%%%%%%%%%%%%%%%%%%%%%%%%%%%%%%%%%%%%%%%%%%%%%%%%%%%%%%%%%%%%%%%%%%%%%%%%%%%%%%%%%%%%%%%%
\subsection{Penalization Scheme for a Generalized BSDE}             \label{nsec3.4}
%%%%%%%%%%%%%%%%%%%%%%%%%%%%%%%%%%%%%%%%%%%%%%%%%%%%%%%%%%%%%%%%%%%%%%%%%%%%%%%%%%%%%%%%%%%%%%%%

To examine the penalization scheme for the generalized BSDE, we define $\barS := S^r\cup\{T\}$ where $S^r=S^r(A)$ is the right support of the driver $A$. Recall that the \emph{right support} of a nondecreasing process $A$ is given by $S^r(A):=\{ t\in \RR_+ : \forall\,\varepsilon>0 \ A_{t+\varepsilon}-A_t>0\}$. Furthermore, we denote by $\bSTF$ the class of all $\FF$-stopping times  $\tau$ taking values in $[0,T]$ and such that $\PP ( \tau \in \barS)=1$. Similarly, $\bSTFt$ is the set of $\FF$-stopping times from $\bSTF$ such that $\PP ( \tau \in \barS \cap [t,T] )=1$.
%\bhyp \label{nhyp1.3}
%The set $\barS:=\xS(\newA )\cup\{T\}$ where $\xS(\newA )$ is the support of the random measure generated by $\newA$.
%\ehyp

Before proceeding to the penalization scheme, we present the following auxiliary result, which is an extension of Lemma 6.1 in \cite{EKPPQ1997}.

\bl  \label{lem3.4x}
Let $\eta$ be an $\FF$-optional, bounded and right-continuous process. Then for any stopping time $\nu\in \bSTF$ we have
\begin{align*}
\lim_{\,n\to\infty}\big[ \xi_T \cE_{\nu,T}(-A^n)+(\I_{\rrb \nu,T\rrb}\eta\cE_{\nu,\cdot}(-A^n)\bcd A^n)_T\big]
=\xi_T\I_{\{\nu=T\}}+\eta_{\nu}\I_{\{\nu<T\}}
\end{align*}
where $A^n := nA$.
\el

\begin{proof}
Since $A^n_\nu-A^n_T=0$ on the event $\{\nu=T\}$ we have
\begin{align*}
\lim_{\,n\to\infty} \cE_{\nu,T}(-A^n)=\lim_{\,n\to\infty} e^{A^n_\nu-A^n_T}=\I_{\{\nu=T\}}.
\end{align*}
Since $\nu$ takes values in the right support we have that $A^n_T - \wt\Gamma_\nu > 0$ on the event $\{\nu<T\}$
and thus $\lim_{\,n\to\infty} (A^n_\nu-A^n_T)=-\infty$. We claim that, on the event $\{\nu<T\}$,
\begin{align}
\lim_{\,n\to\infty}\int_{\rrb\nu,T\rrb}\eta_s e^{A^n_{\nu}-A^n_s}\,d\wtGamma^n_s
=\eta_{\nu}\I_{\{\nu<T\}}\label{dirac}
\end{align}
since the sequence of bounded and positive measures $\mu_n := \I_{\rrb\nu,T\rrb} e^{A^n_{\nu}-A^n_s}\,dA^n_s,\, n \in \NN$ converges to the Dirac measure at $\nu$ on the event $\{\nu<T\}$, that is, to the measure $\mu:=\I_{\{\nu<T\}}\delta_{\nu}$.
Equality \eqref{dirac} can be formally established using the time change on $[0,T]$ generated by $A$. To this end, we define the increasing, right-continuous process $C$ by $C_s=\inf\{t\in\RR_+:A_t>s\}$. Then an application of the time change formula (see, e.g., Chapter 0 in Revuz and Yor \cite{RY1999}) gives
\begin{align*}
\int_{\rrb\nu,T\rrb}\eta_se^{A^n_{\nu}-A^n_s}\,dA^n_s
&=\int_{\rrb\nu,T\rrb}\eta_s ne^{n(A_{\nu}-A_s)}\,dA_s\\
&=\int_0^\infty \I_{\{\nu< C_s \leq T\}}  \eta_{C_s} n e^{n(A_{\nu}-s)}\,ds
\end{align*}
where to obtain the second equality, we have used the equality $A_{C_s} =s$, which holds since $A$ is a continuous process. From the fact that $\{\nu < C_s\} \subseteq \{A_\nu \leq s\}$ and the change of variable $u = s- A_\nu$ we obtain
\begin{align*}
\int_0^\infty\I_{\{\nu< C_s \leq T\}}\eta_{C_s} ne^{n(A_{\nu}-s)}\,ds &=\int^\infty_0\I_{\{\nu< C_{A_\nu+ u}\leq T\}}\eta_{C_{A_\nu+ u}} n e^{-nu}\,du\\
&=\EE_X\Big[\I_{\{\nu< C_{A_{\nu+ X/n}}\leq T\}}\eta_{C_{A_{\nu+ X/n}}}\Big]
\end{align*}
where in the last equality we have used the observation that $ne^{-nu},\, u>0$ is the density of $\frac{1}{n}X$ where $X \sim \exp(1)$ and is independent of $\cG_\infty$.

Since on the event $\{\nu< T\}$, the stopping time $\nu$ takes values in the right support of $A$ and $\nu < \nu + X/n$, we deduce that $A_{\nu + X/n} \geq A_{\nu + X/(n+1)}  > A_{\nu}$ for a sufficiently large $n \in \mathbb{N}$. This observation, together with the right-continuity of the processes $\eta$, $C$ and $A$, allows us to conclude that
\begin{align*}
\lim_{\,n\rightarrow \infty} \EE_X\Big[\I_{\{\nu< C_{A_\nu+ X/n} \leq T\}}\eta_{C_{A_\nu+ X/n}}\Big]=\I_{\{\nu \leq C_{A_\nu}<T\}}\eta_{C_{A_\nu}}.
\end{align*}
Finally, since on $\{\nu<T\}$ the stopping time $\nu$ takes values in the right support of $A$ we have $C_{A_\nu}=\inf\{s: A_s>A_\nu\}=\nu$, which allows us to conclude that $\I_{\{\nu \leq C_{A_\nu}<T\}}\eta_{C_{A_\nu}}=\I_{\{\nu<T\}}\eta_{\nu}$.
\end{proof}

\bt \label{th3.1}
Let Assumptions \ref{nass1.1} and \ref{nass1.2} hold, $\xi_T$ be bounded, and $\eta$ be an $\FF$-optional, bounded, nonnegative, and right-continuous process. Consider the sequence of unique solutions $(Y^n,Z^n)\in \cS^2(\FF) \times \cH^2(M)$ to the generalized BSDE
\begin{align} \label{eq3.6}
Y^n_t=\xi_T-\int_{]t,T]} Z^n_s\,dM_s+\int_{]t,T]} n(\eta_s-Y^n_s)^+\,dA_s.
\end{align}
Then the sequence $Y^n$ is increasing and converges almost surely to the process $V$, which satisfies, for every $t\in [0,T]$,
\begin{equation} \label{eq3.7}
V_t=\esssup_{\tau \in\bSTFt}\EE\big[\xi_{T}\I_{\{\tau=T\}}+\eta_{\tau}\I_{\{\tau< T\}}\,|\,\cF_t\big].
\end{equation}
\et

\begin{proof}
It is known from Proposition \ref{npro1.4} that the generalized BSDE \eqref{eq3.6} has a unique solution $(Y^n,Z^n)$ in $\cS^2(\FF) \times \cH^2(M)$.
We note that the sequence $Y^n$ of processes is monotonically increasing as $n\rightarrow \infty$ (see Corollary \ref{ncor1.2}) and the limit $Y=\lim_{\,n\rightarrow \infty}Y^n$ is well defined.

\noindent {\it Step 1.}  Our first goal is to show that the c\`adl\`ag, $\FF$-adapted process $Y^n$  satisfies, for every $n \in \NN$ and all $t\in [0,T]$,
\begin{align} \label{eq3.8}
Y_t^n=\esssup_{\tau\in \bSTFt}\EE\big[\xi_{T}\I_{\{\tau=T\}}+(\eta_{\tau}\land Y_{\tau}^n)
\I_{\{\tau< T\}}\,|\,\cF_t\big].
\end{align}
For a fixed $n\in \NN$, we set $\gamma^n_t :=  \xi_{T}\I_{\{t=T\}}+(\eta_{t}\land Y_{t}^n)\I_{\{ t < T\}}$ and we observe that
$$\gamma^n_t = (\xi_{T}\wedge Y_T^n)\I_{\{t=T\}}+(\eta_{t}\land Y_{t}^n)\I_{\{ t < T\}} \leq Y_t^n.$$
Furthermore, the GBSDE \eqref{eq3.6} can be represented as
\begin{align} \label{eq3.9}
Y^n_t=\xi_T-\int_{]t,T]}Z^n_s\,dM_s+K^n_T-K^n_t
\end{align}
where the $\FF$-adapted, continuous, nondecreasing process $K^n$ is given by
\begin{align*}
K^n_t:=\int_{]0,t]} n(\eta_s-Y^n_{s})^+\,dA_s .
\end{align*}
Recall that we assumed that $\eta$ (and hence also $\gamma^n$) is a right-continuous process. We claim that \eqref{eq3.8} is valid,
that is, for every $t\in [0,T]$,
\begin{align} \label{eq3.10}
Y_t^n=\esssup_{\tau\in \bSTFt}\EE\big[\gamma^n_{\tau}\,|\,\cF_t\big].
\end{align}
Equality \eqref{eq3.10} is obvious for $t=T$ so it suffices to consider any $t<T$. We have, for any $\tau \in \bSTFt$,
\begin{align*}
Y_t^n =\EE\big[ Y^n_{\tau}+K^n_{\tau}-K^n_t  \,|\,\cF_t\big]
\geq \EE\big[ Y^n_{\tau} \,|\,\cF_t\big] \geq \EE\big[ \gamma^n_{\tau} \,|\,\cF_t\big],
\end{align*}
which implies that $Y_t^n \geq \esssup_{\tau\in \bSTFt}\EE\big[\gamma^n_{\tau}\,|\,\cF_t\big]$.

For the converse inequality, we define the stopping time $\tau_t:=\inf\{s\in [t,T]\,|\, K^n_s-K^n_t>0\}$, which belongs to $\bSTFt$,
and we observe that $K^n_{\tau_t}-K^n_t=0$ due to the continuity of $K^n$.
Furthermore, on the event $\{\tau_t< T\}$ we have $K^n_s>K^n_{\tau_t}$ on $\rrb \tau_t,T\rrb$, which
entails the inequality $\limsup_{u \downarrow \tau_t}(\eta_u-Y^n_u)\geq 0$ and thus, since the processes
$\eta$ and $Y^n$ are right-continuous, we conclude that $\eta_{\tau_t}-Y^n_{\tau_t}\geq 0$,
which in turn implies that $Y^n_{\tau_t}= \gamma^n_{\tau_t}$.
On the event $\{\tau_t = T\}$ we have that $K^n_{\tau_t}-K^n_t=K^n_T-K^n_t=0$ and
$Y^n_{\tau_t}=\gamma^n_{\tau_t}$. From \eqref{eq3.9}, we now get
\begin{align*}
Y^n_t=Y^n_{\tau_t}-\int_{\rrb t,\tau_t\rrb}Z^n_s\,dM_s=\gamma^n_{\tau_t}-\int_{\rrb t,\tau_t\rrb}Z^n_s\,dM_s
= \EE\big[ \gamma^n_{\tau_t}\,|\,\cF_t\big],
\end{align*}
which leads to the inequality $Y_t^n \leq \esssup_{\tau\in \bSTFt}\EE\big[\gamma^n_{\tau}\,|\,\cF_t\big]$ since $\tau_t \in \bSTFt$.
We conclude that equality \eqref{eq3.10} (or, equivalently, \eqref{eq3.8}) holds.

\noindent {\it Step 2.}  We set $\gamma_t:=\xi_{T}\I_{\{t=T\}}+\eta_{t}\I_{\{t< T\}}$ and we show that $\lim_{\,n\rightarrow\infty} Y_t^n = V_t$ where
\begin{align}  \label{eq3.11}
V_t :=\esssup_{\nu \in \bSTFt}  \EE\big[\gamma_{\nu}\,|\,\cF_t\big].
\end{align}
It follows from  \eqref{eq3.8} and \eqref{eq3.11} that $Y^n \leq V$ and  $Y^n \uparrow Y \leq V$ where the $\FF$-optional process $Y$ is given by $Y=\lim_{\,n\rightarrow +\infty} Y^n$. Furthermore, it is clear that the process $Y$ is nonnegative and belongs to $\cS^2(\FF)$ since it is dominated by
the process $V$ belonging to $\cS^2(\FF)$.

It is clear that $Y_T =\xi_T =  \gamma_T$ and thus it suffices to show that $Y_\nu\I_{\{\nu < T\}} \geq \eta_\nu\I_{\{\nu < T\}}$ for any
$\FF$-stopping time $\nu$ taking values in $\barS$. From the monotone convergence theorem and the comparison property established
in Proposition \ref{npro1.1} we obtain, for every $0\leq \tau \leq \nu \leq T$,
\begin{gather*}
 \EE\big[ Y_\nu \,|\, \cF_\tau \big] = \lim_{\,n\rightarrow \infty} \EE\big[Y^n_\nu \,|\,\cF_\tau \big]
 \leq \lim_{\,n\rightarrow \infty} Y^n_\tau = Y_\tau.
\end{gather*}
Using the fact that $V$ is bounded, we deduce that $Y$ is a bounded strong supermartingale and, as a consequence of Theorem 2 in Mertens \cite{M1972} (or, more specifically, the Lemma on page 51 of \cite{M1972}), we have $Y \geq Y_{+}$. Next, from the form of $Y^n$ we have,
for every $\tau \in \STF$,
\begin{gather*}
\frac{1}{n}\EE[Y^n_\tau]=\frac{1}{n}\EE[\xi_T]+\EE \Big[\int_{\rrb\tau,T\rrb} (\eta_s-Y^n_s)^+\,dA_s\Big].
\end{gather*}
By letting $n$ go to $\infty$ and applying the dominated convergence theorem, we obtain the equality $\int_{\rrb\tau,T\rrb}(\eta_s-Y_s)^+\,dA_s=0$. We claim that the last equality implies that $Y_\nu\I_{\{\nu<T\}} \geq \eta_\nu\I_{\{\nu<T\}}$ for all $\nu\in \bSTF$. Suppose, on the contrary, that there exists a stopping time $\nu\in \bSTF$ such that the event $E = \{\eta_\nu > Y_\nu \geq Y_{\nu+}\}\cap \{\nu < T\}$ has a positive probability. Then we deduce
from the right-continuity of $\eta$ that there exists $\varepsilon$, which may depend on $\omega \in E$, such that $\eta_{t} - Y_{t} > 0$ for all $t \in \llb \nu, \nu + \varepsilon \rrb  \cap \llb \nu, T\rrb$. However, since $\nu \in \bSTF$, we deduce that $A_{\nu + \varepsilon} - A_\nu > 0$ for almost every $\omega\in E$ and this contradicts the equality $\int_{\rrb\nu,T\rrb} (\eta_s-Y_s)^+\,dA_s= 0$. We this conclude that $Y$ dominates $\gamma$ on $\barS$.
%
%Furthermore, to show that $Y$ dominates $\gamma$ on $\barS$ we note that, for every $\tau \in \bSTF$,
%\begin{align*}
%Y_\tau^n& \geq  \barY^n_\tau = \xi_T+ \int_{\rrb \tau, T\rrb} n(\eta_s- \barY^n_s)\,dA_s-\int_{\rrb\tau,T\rrb} \barZ^n_s\,dM_s
%\end{align*}
%
%
%where the right-hand side is a linear BSDE and we have, for every $\tau \in \bSTFt$,
%\begin{align*}
%\barY_\tau^n= \EE\big[\xi_{T}\cE_{\tau,T}(-A^n)+(\I_{\rrb \tau,T\rrb}\eta\cE_{\tau,\cdot}(-A^n)\bcd A^n)_{T}\,|\,\cF_\tau\big]
%\end{align*}
%where $A^n= n A$. We observe that $\cE_{\tau,T}(-A^n)$ converges to $\I_{\{\tau=T\}}$ as $n\rightarrow \infty$ and thus, using Lemma \ref{lem3.4x} and the right-continuity of $\eta$, we obtain
%\begin{align}
%\barY_\tau=\lim_{\,n\rightarrow \infty}Y^n_{\tau} \geq \lim_{\,n \rightarrow \infty} \barY^n_{\tau}=\eta_{\tau}\I_{\{\tau= \sigma\}}+\eta_{\tau}\I_{\{\tau< \sigma\}} = \gamma_\tau. \label{VR}
%\end{align}
%Thus $\barY$ dominates $\gamma$ on $\barS(A)$. To conclude that $\barY = V$, it now suffices to show that  $\barY \geq V$.

Finally, by using the fact that the essential supremum and conditional expectation can be interchanged, we can conclude that $V$ is the smallest strong supermartingale dominating $\gamma$ on $\barS(A)$,
%From the monotone convergence theorem we obtain, for every $0\leq \tau \leq \nu \leq T$,
%\begin{gather*}
% \EE\big[ Y_\nu \,|\, \cF_\tau \big] = \lim_{\,n\rightarrow \infty} \EE\big[Y^n_\nu \,|\,\cF_\tau \big]
% \leq \lim_{\,n\rightarrow \infty} Y^n_\tau = Y_\tau.
%\end{gather*}
which in turn implies that,
for every $t\in[0,T]$,
\begin{gather*}
Y_t \geq \esssup_{\nu \in \bSTFt} \EE\big[\gamma_{\nu}\,|\,\cF_t\big]=V_t ,
\end{gather*}
as was required to show.
\end{proof}

%%%%%%%%%%%%%%%%%%%%%%%%%%%%%%%%%%%%%%%%%%%%%%%%%%%%%%%%%%%%%%%%%%%%%%%%%%%%%%%%%%%%%%%%%%%%%
%%%%%%%%%%%%%%%%%%%%%%%%%%%%%%%%%%%%%%%%%%%%%%%%%%%%%%%%%%%%%%%%%%%%%%%%%%%%%%%%%%%%%%%%%%%%%
\section{Reflected Generalized BSDEs}             \label{nsec3}
%%%%%%%%%%%%%%%%%%%%%%%%%%%%%%%%%%%%%%%%%%%%%%%%%%%%%%%%%%%%%%%%%%%%%%%%%%%%%%%%%%%%%%%%%%%%%
%%%%%%%%%%%%%%%%%%%%%%%%%%%%%%%%%%%%%%%%%%%%%%%%%%%%%%%%%%%%%%%%%%%%%%%%%%%%%%%%%%%%%%%%%%%%%

Our next goal is to study reflected generalized BSDEs (reflected GBSDEs or, simply, RGBSDEs) with a lower obstacle given by a predetermined
bounded and $\FF$-optional process $\zeta$. It is known that if the process $\zeta$ is bounded and $\FF$-optional, then the process $\bzeta$ given by $\bzeta_t:=\limsup_{s \uparrow t,\,s<t}\zeta_t$ for all $t \in ]0,T]$ is known to be $\FF$-predictable and left-upper-semicontinuous
(see, e.g., \cite{GIOQ2020} or Theorem 90 on page 225 in \cite{DM1975}). The process $\bzeta$ is called the {\it left-upper-semicontinuous envelope} of $\zeta$.

Recall that by $\cK$ (resp., $\bcK$) we denote the class of all c\`adl\`ag, nondecreasing, $\FF$-predictable (resp., l\`adl\`ag, nondecreasing, $\FF$-predictable) processes and $\STF$ (resp., $\STFp$) stands for the class of all $\FF$-stopping times (resp., $\FF$-predictable stopping times) $\tau$ taking values in $[0,T]$. The following definition is consistent with the classical case where $M=W$ is a Brownian motion and the
driver $A_t=t$ for every $t\in [0,T]$.

\bd \label{ndef1.4}
A triplet $(Y,Z,K)$ is a {\it solution to the reflected generalized BSDE} with data
$(M,\newA,\gen,\zeta)$ if $Y$ is an $\FF$-adapted, l\`adl\`ag process, $Z$ is an $\FF$-predictable process,
$K$ is a nondecreasing, $\FF$-predictable, l\`adl\`ag  process, and the following conditions are met
\begin{equation} \label{RGBSDE}
\left\{ \begin{array} [c]{ll}
Y_\tau=\zeta_T-\int_{\rrb \tau ,T\rrb}Z_s\,dM_s+\int_{\rrb\tau,T\rrb} \gen (s,Y_s)\,d\newA_s+ K_T- K_\tau, \ \forall \, \tau \in \STF , \medskip\\
Y_t \geq \zeta_t , \ \forall\, t \in [0,T], \quad \big(\I_{\{Y_{-}>\bar{\zeta} \}}\bcd K^c\big)_T=0 ,  \medskip\\
(Y_{\tau-}-\bar{\zeta}_{\tau})\Delta K^{d}_{\tau}=0,\ \forall\, \tau \in \STFp,\quad
(Y_{\tau}-\zeta_{\tau})\Delta^+ K^g_{\tau}=0, \ \forall\, \tau \in \STF,
\end{array} \right.
\end{equation}
where the It\^o integral $\int_{]0,t]}Z_s\,dM_s$ is a martingale and the integral $\int_{]0,t]} \gen (s,Y_s)\,d\newA_s$ is
an $\FF$-adapted, continuous process of finite variation.
\ed

Our first goal is to establish in Section \ref{nsec3.1} two variants of the comparison theorem for a reflected GBSDE and obtain in Section \ref{nsec3.2}
some useful {\it a priori} estimates. Subsequently, we will study in Section \ref{nsec3.3} the existence and uniqueness of a solution $(Y,Z,K)$ to the RGBSDE \eqref{RGBSDE} in the space $\cS^2(\FF) \times \cH^2(M)\times \bcK$. Finally, for a given $\FF$-optional and bounded process $\eta$ we will consider in Section \ref{nsec3.4a} the sequences of
penalized RGBSDEs
\begin{align*} % \label{veq6.9v}
Y^n_{\tau}= \zeta_T-\int_{\rrb\tau,T\rrb} Z^n_s\,dM_s+\int_{\rrb\tau,T\rrb}  n(Y^n_s-\eta_s)^+\,d\newA_s+K^n_T-K^n_{\tau}
\end{align*}
and
\begin{align*}  % \label{beq6.9v}
\wtY^n_{\tau}= \zeta_T-\int_{\rrb\tau,T\rrb}  \wtZ^n_s\,dM_s - \int_{\rrb\tau,T\rrb}  n(\wtY^n_s-\eta_s)^+\,d\newA_s+\wtK^n_T-\wtK^n_{\tau}
\end{align*}
with the lower obstacle $\zeta$ and we will examine the limits $Y:=\lim_{\,n\to\infty}Y^n$ and $\wtY:=\lim_{\,n\to\infty}\wtY^n$.
It will be shown that the process $Y$ (resp. the process $\wtY$) can be interpreted as the value process of a particular optimal
stopping problem (resp., a Dynkin game) where the right support of the driver $A$ plays an important role in the specification
of the respective reward process.

%%%%%%%%%%%%%%%%%%%%%%%%%%%%%%%%%%%%%%%%%%%%%%%%%%%%%%%%%%%%%%%%%%%%%%%%%%%%%%%%%%%%%%%%%%%%%%
\subsection{Comparison Theorems for a Reflected GBSDE}       \label{nsec3.1}
%%%%%%%%%%%%%%%%%%%%%%%%%%%%%%%%%%%%%%%%%%%%%%%%%%%%%%%%%%%%%%%%%%%%%%%%%%%%%%%%%%%%%%%%%%%%%%

We aim to show that the comparison property of Proposition \ref{npro1.1} can be extended from a GBSDE to the case of a reflected GBSDE.
The first variant of the comparison theorem for a reflected GBSDE is established under an additional assumption
that the filtration $\FF$ is quasi-left-continuous so that the martingale $M$ does not jump
at any $\FF$-predictable stopping time.

\bp \label{npro1.5}
Suppose that the filtration $\FF$ is quasi-left-continuous. Let the mappings $\gen ,h:\Omega \times [0,T]\times \RR \to \RR$ be such that the reflected GBSDEs with the lower obstacles $\zeta$ and $\wtzeta$
\begin{align*}
Y_{\tau}=\zeta_T-\int_{\rrb \tau ,T\rrb}Z_s\,dM_s+\int_{\rrb \tau ,T\rrb} \gen (s,Y_s)\,d\newA_s + K_T-K_{\tau}
\end{align*}
and
\begin{align*}
\wtY_{\tau}=\wtzeta_T-\int_{\rrb \tau ,T\rrb}\wtZ_s\,dM_s+\int_{\rrb \tau ,T\rrb} h(s,\wtY_s)\,d\newA_s + \wtK_T - \wtK_{\tau}
\end{align*}
have % unique
solutions $(Y,Z,K)$ and $(\wtY,\wtZ,\wtK)$, respectively. % in $\cS^2(\FF)\times \cH^2(M) \times \bcK$.
If $\zeta \ge \wtzeta $ and the functions $\gen$ and $h$ satisfy the following conditions: \\
(i) $\gen (\omega,t,y)\ge h(\omega,t,y)$ for every $(\omega,t,y) \in \Omega \times [0,T] \times \RR$, \\
(ii) $\gen (\omega,t,\cdot )$ % or $h(\omega,t,\cdot)$
is a nonincreasing function for every $(\omega,t) \in \Omega \times [0,T],$ \\
then the inequality $Y\ge \wtY $ is valid.
\ep

\begin{proof}
We will extend the proof of Proposition \ref{npro1.1}. It is important to observe that the processes $Y$ and $\wtY$ introduced in the
statement of Proposition \ref{npro1.5} are l\`adl\`ag whereas they were c\`adl\`ag in Proposition \ref{npro1.1},
For a fixed $\varepsilon>0$, we define the $\FF$-stopping time $\tau^{\varepsilon}:=\inf\,\{t\ge 0 : Y_t \leq \wtY_t - \varepsilon\}$ where, by convention, $\inf \emptyset = T$.  As in the proof of Proposition \ref{npro1.1}, we argue that if for all $\varepsilon>0$ we have $\PP (\tau^{\varepsilon}=T)=1$, then the asserted inequality $Y\ge \wtY$ holds.

Let us now assume that the inequality $Y\geq \wtY$ does not hold. Then there exists $\varepsilon >0$ such that $\PP (\tau^{\varepsilon} <T)>0$. We
fix $\varepsilon$ and we denote $E:=\{\tau^{\varepsilon} <T\} \in \cF_{\tau^{\varepsilon}}$. Next, we define $\tau:=\tau^{\varepsilon}\I_E + T\I_{E^c}$
so that $\{\tau<T\}=\{\tau^{\varepsilon} <T\}\in \cF_{\tau}$.

\vskip 5 pt
\noindent {\it Step 1.} We start by showing that the inequality $Y_{\tau+}<\wtY_{\tau+}$ holds on $E$.
To this end, let us consider any event $\omega$ from $E$. We have that either (a.1) $Y_{\tau} \leq \wtY_{\tau}-\varepsilon$ so that
$Y_{\tau}<\wtY_{\tau}$ or (a.2) $Y_{\tau}> \wtY_{\tau}-\varepsilon$ but $Y_{\tau+} \leq \wtY_{\tau+}-\varepsilon$ so that $Y_{\tau+}<\wtY_{\tau+}$, which is the desired inequality.

Hence it remains to show that in case (a.1) we also have that $Y_{\tau+}<\wtY_{\tau+}$. To this end,  we first observe that
$\wtY_{\tau} > \wtY_{\tau}-\varepsilon > Y_{\tau} \geq \zeta_{\tau} \geq \wtzeta_{\tau}$ and thus the process $\wtK$ is right-continuous at $\tau$ (from the respective Skorokhod condition), which in turn implies that $\wtY$ is right-continuous at $\tau$. Recall that $Y_{\tau}=\zeta_{\tau}\vee Y_{\tau+}$. If $Y_{\tau}>\zeta_{\tau}$, then $Y$ is also right-continuous at $\tau$ and thus $\wtY_{\tau+}=\wtY_{\tau}> Y_{\tau}=Y_{\tau+}$. Finally, if $Y_{\tau}=\zeta_{\tau}$, then $Y_{\tau+}\leq \zeta_{\tau}$ and thus $\wtY_{\tau+}=\wtY_{\tau}>\zeta_{\tau}\geq Y_{\tau+}$.

We have thus shown that the inequality $Y_{\tau+}>\wtY_{\tau+}$ holds on $E$, for almost all $\omega$. Then we define the $\FF$-stopping
time $\nu:=\inf\,\{t\ge \tau :Y_t\geq\wtY_t\}$ and we note that $\nu\leq T$ since, by assumption, $Y_T =\zeta_T \geq \wtzeta_T =\wtY_T$. Since $Y_{\tau+}>\wtY_{\tau+}$ on $E$ it is clear that the interval $\rrb \tau,\nu\rrb $ is nonempty on the event $E=\{\tau <T\} \in \cF_{\tau}$.
It is also worth noting that $E$ belongs also to $\cF_{\nu-}$ since $\tau < \nu$ on $E$ (see Proposition 2.4 in Nikeghbali \cite{N2006}).

\vskip 5 pt
\noindent {\it Step 2.} Our next goal is to show that the inequality $Y_{\nu} \geq \wtY_{\nu}$ is satisfied.
It manifestly holds on the event $\nu =T$ and thus it suffices to shown that it is valid on $\{\nu < T\}$ as well.
If $Y_{\nu} \geq \wtY_{\nu}$, then the desired inequality manifestly holds and thus it suffices to examine the event $\{Y_{\nu}<\wtY_{\nu},\, Y_{\nu+} \geq \wtY_{\nu+}\}$. We will show by contradiction that the probability of that event is null.
Since $Y_{\nu} = \zeta_{\nu}\vee Y_{\nu+}$ it suffices to consider two cases: (b.1) $Y_{\nu}>Y_{\nu+}$ and (b.2) $Y_{\nu}= Y_{\nu+}$.

In case (b.1), we have $Y_{\nu} = \zeta_{\nu}$ and thus $\wtY_{\nu}> Y_{\nu}= \zeta_{\nu} \geq \wtzeta_{\nu}$. This implies that $\wtY$ is right-continuous at $\nu$, which in turn yields $\wtY_{\nu}=\wtY_{\nu+} > Y_{\nu}>Y_{\nu+}$ and hence contradicts the assumption that $Y_{\nu}<\wtY_{\nu}$.

In case (b.2), we have $Y_{\nu} = Y_{\nu+}$ and thus $\wtY_{\nu}> Y_{\nu}=Y_{\nu+} \geq \wtY_{\nu+}$,
which implies $\wtY_{\nu}=\wtzeta_{\nu}> Y_{\nu} \geq  \zeta_{\nu}$. Hence $\wtzeta_{\nu}>\zeta_{\nu}$, which is a contradiction since,
by assumption, the inequality $\wtzeta \leq \zeta$ holds. We thus see that the inequality $Y_{\nu}\geq \wtY_{\nu}$ is proven.

\vskip 5 pt
\noindent {\it Step 3.} We are now ready to show that if $\PP(E)>0$ then a contradiction arises.
From the first step, we deduce that there exists a sufficiently small constant $\delta = \delta (\varepsilon)>0$ such that
$\PP (C):=\PP ( \tau +\delta < \nu)>0$ and $Y_{\tau+\delta}<\wtY_{\tau+\delta}$. Hence we define the $\FF$-stopping time $\sigma:=(\tau+\delta )\I_E + T\I_{E^c}$ and we consider the interval $\llb \sigma,\nu \rrb$. We henceforth work on the event $E \in \cF_{\sigma}$.
As in Section \ref{nsec2.1} we denote
\begin{align*}
\Gen_t(Y):=\int_{]0,t]} \gen (s,Y_s)\,d\newA_s ,\quad H_t(\wtY):=\int_{]0,t]} h(s,\wtY_s)\,d\newA_s
\end{align*}
and we also write $U:=G(Y)$ and $\wtU:=H(\wtY)$ so that the process $\ov{U}:=U-\wtU$ satisfies
\begin{align*}
\ov{U}_t=\Gen_t(Y)-H_t(\wtY)=(\Gen_t(Y)-\Gen_t(\wtY))+(\Gen_t(\wtY)-H_t(\wtY)).
\end{align*}
Since $\gen (\omega,t,y) \ge h(\omega,t,y)$ for all $(\omega ,t,y) \in \Omega \times [0,T]\times \RR$, the inequality $Y <\wtY$ holds on  $\llb\sigma,\nu \llb$ and for every $(\omega ,t) \in \Omega \times [0,T]$ the function $\gen (\omega,t,\cdot )$ is nonincreasing,
it is easy to check that $\I_E \ov{U}$ is a continuous and nondecreasing process on $\llb\sigma,\nu\rrb$.
We observe that on $\llb\sigma,\nu \rrb$
\begin{align*}
Y_t-\wtY_t =Y_{\nu}-\wtY_{\nu}-\int_{\rrb t,\nu\rrb}(Z_s-\wtZ_s)\,dM_s+(K_{\nu}-K_t)-(\wtK_{\nu}-\wtK_{\nu-})+\ov{U}_{\nu}-\ov{U}_t
\end{align*}
where in fact the process $\wtK$ is constant on $\llb\sigma,\nu\llb$ since $\wtY > Y \geq \zeta \geq \wtzeta$ on $\llb\sigma,\nu\llb$
and thus $\wtK_{t} = \wtK_{\nu-}$ for $t \in \llb\sigma,\nu\llb$. Consequently,
\begin{align}
Y_t-\wtY_t =Y_{\nu}-\wtY_{\nu} + \Delta K_{\nu} - \Delta \wtK_{\nu} -\int_{\rrb t,\nu \rrb}(Z_s-\wtZ_s) \,dM_s+(K_{\nu-}-K_t) +\ov{U}_{\nu}-\ov{U}_t .
\end{align}
By taking the $\cF_\sigma$-conditional expectation, we obtain
\begin{align*}
Y_\sigma-\wtY_\sigma  & \geq \EE\big[Y_{\nu}-\wtY_{\nu} + \Delta K_{\nu} - \Delta \wtK_{\nu} \, |\, \cF_\sigma\big] \\
					  & = \EE\big[(Y_{\nu}-\wtY_{\nu} + \Delta K_{\nu})\I_{\{\Delta \wtK_\nu = 0 \}} \, |\, \cF_\sigma\big]
                         + \EE\big[(Y_{\nu}-\wtY_{\nu} + \Delta K_{\nu} - \Delta \wtK_{\nu})\I_{\{\Delta \wtK_\nu > 0 \}} \, |\, \cF_\sigma\big] \\
					  & \geq \EE\big[(Y_{\nu}-\wtY_{\nu} + \Delta K_{\nu} - \Delta \wtK_{\nu})\I_{\{\Delta \wtK_\nu > 0 \}} \, |\, \cF_\sigma\big]
\end{align*}
where in the last inequality we have used the facts that $Y_{\nu}- \wtY_{\nu}\geq 0$ and $\Delta K_{\nu} \geq 0$. Next, we show that
\begin{align*}
\EE\big[(Y_{\nu}-\wtY_{\nu} + \Delta K_{\nu} - \Delta \wtK_{\nu})\I_{\{\Delta \wtK_\nu > 0 \}} \, |\, \cF_\sigma\big] = 0.
\end{align*}
We first notice that on the event $\{\Delta \wtK_\nu > 0 \}$
\begin{align*}
(Y_{\nu}-\wtY_{\nu} + \Delta K_{\nu} - \Delta \wtK_{\nu})\I_{\{\Delta \wtK_\nu > 0 \}}= (Z_{\nu}-\wtZ_{\nu})\,\Delta M_{\nu}\I_{\{\Delta \wtK_\nu > 0 \}}
\end{align*}
since it is easy to check that $Y_{\nu-}-\wtY_{\nu-}=0$ on $\{\Delta \wtK_{\nu}> 0\}$. The filtration $\FF$ is assumed to be quasi-left-continuous and thus the equality $\Delta M_{\tau} = 0$ holds for any $\FF$-predictable stopping time $\tau$ (see Theorem 5.39 in \cite{HWY1992}). Since the process $\wtK$ is strongly $\FF$-predictable, by Theorem 3.33 in \cite{HWY1992}, the set $\{\Delta \wtK > 0\} = \{\Delta \wtK^d > 0\}$ is included in the union of the graphs of a family of $\FF$-predictable stopping times. Hence we can conclude that $\Delta M\I_{\{\Delta \wtK > 0\}} = 0$ and thus $Y_\sigma \geq \wtY_\sigma$. However, $Y_{\sigma}-\wtY_{\sigma}<0$ on $E\in \cF_{\sigma}$ and thus $Y_{\sigma}=\wtY_{\sigma}$ on $E$, which contradicts the definition of $\sigma$.
\end{proof}

In the second variant of the comparison theorem we relax the assumption that the filtration $\FF$ is quasi-left-continuous.
However, the method of proof requires to postulate that the random variable $A_T$ is bounded, and not merely square-integrable.
Finally, we denote the stochastic exponential of a semimartingale $X$ by $\cE(X)$  (see Theorem 5.1 in \cite{G1985})  and we set $\cE_{s,t}(X)=\cE_t(X)/\cE_s(X)$ for every $s\leq t$.

\bp \label{npro1.5x}
Suppose that the process $A$ is bounded. Let the mappings $\gen ,h:\Omega \times [0,T]\times \RR \to \RR$ be such that the reflected GBSDEs with the lower obstacles $\zeta$ and $\wtzeta$
\begin{align*}
Y_{\tau}=\zeta_T-\int_{\rrb \tau ,T\rrb}Z_s\,dM_s+\int_{\rrb \tau ,T\rrb} \gen (s,Y_s)\,d\newA_s+ K_T-K_{\tau}
\end{align*}
and
\begin{align*}
\wtY_{\tau}=\wtzeta_T-\int_{\rrb \tau ,T\rrb}\wtZ_s\,dM_s+\int_{\rrb \tau ,T\rrb} h(s,\wtY_s)\,d\newA_s + \wtK_T - \wtK_{\tau}
\end{align*}
have unique solutions $(Y,Z,K)$ and $(\wtY,\wtZ,\wtK)$ in $\cS^2(\FF)\times \cH^2(M) \times \bcK$.
Suppose that the functions $\gen$ and $h$ satisfy the following conditions: \\
(i) $\gen (\omega,t,y)\ge h(\omega,t,y)$ for every $(\omega,t,y) \in \Omega \times [0,T] \times \RR$, \\
(ii) $\gen (\omega,t,\cdot )$ % or $h(\omega,t,\cdot)$
is a nonincreasing function for every $(\omega,t) \in \Omega \times [0,T].$ \\
If the obstacles satisfy $\zeta \ge \wtzeta $ and are right-upper-semicontinuous, then the inequality $Y\ge \wtY $ is valid.
\ep

\begin{proof}
{\it Step 1.} We first show that the process $Y$ given by
\begin{align*}
Y_{\tau}=\zeta_T-\int_{\rrb\tau ,T\rrb}Z_s\,dM_s+\int_{\rrb \tau ,T\rrb} g(s,Y_s)\,dA_s + K_T-K_{\tau}
\end{align*}
is a strong $\cE^g$-supermartingale where the nonlinear evaluation $\cE^g$ (see, e.g., Peng \cite{P2004a,P2004})
is defined through the solution to the GBSDE
\begin{align}
Y'_{\tau}=\zeta_T-\int_{\rrb\tau ,T\rrb} Z'_s\,dM_s+\int_{\rrb \tau ,T\rrb} g(s,Y'_s)\,dA_s . \label{yprime}
\end{align}
We fix $\sigma\in\STF$ and we denote $\barY_s:=Y_s-\cE^g_{s,\sigma}(Y_\sigma)$, $\barZ:=Z-Z'$ and $\overline{g}:=g(\cdot,Y)-g(\cdot,Y')$. We define $\rho := (\overline{g}/\barY)\I_{\{\barY \neq 0\}}$ and we note that the process $\rho$ is clearly bounded since $g$ is a Lipschitz continuous function in $y$ so that the continuous process $X_t:=\int_{]0,t]} \rho_s\,dA_s$ is well defined.  The 
%Gal'\v{c}uk-Lenglart 
integration by parts formula yields, for every $\FF$-stopping times $\tau \leq \sigma $,
\begin{align*}
\mathcal{E}_{\tau,\tau}(X)\barY_\tau
				& = \int_{\rrb \tau ,\sigma\rrb} \mathcal{E}_{\tau,s-}(X)\overline{g}_s\, dA_s - \int_{\rrb \tau ,\sigma\rrb} \mathcal{E}_{\tau,s}(X)\rho_s\barY_s \, dA_s \\
				&\quad  + \int_{\rrb \tau ,\sigma\rrb}  \mathcal{E}_{\tau,s}(X)\, dK^r_s + \int_{\llb \tau ,\sigma\llb}  \mathcal{E}_{\tau,s}(X)\, dK^g_{s+} \\
				 & \quad - \int_{\rrb \tau ,\sigma\rrb} \mathcal{E}_{\tau,s}(X)\barZ_s\, dM_s.
\end{align*}
Since the processes $\rho$ and $A$ are bounded the stochastic exponential $\cE(X)$ belongs to $\mathcal{S}^2(\FF)$
and thus, by the Kunita-Watanabe inequality, the stochastic integral with respect to $M$ is a uniformly integrable martingale. By taking the $\cF_\tau$-conditional expectation and using the assumption that $\overline{g}_\tau =\rho_\tau\barY_\tau$, we obtain
\begin{align*}
\barY_\tau=\cE_{\tau,\tau}(X)\barY_\tau=\EE\Big[\int_{\rrb\tau ,T\rrb}\cE_{\tau,s-}(X)\,dK^r_s
+ \int_{\llb \tau ,T\llb}\cE_{\tau,s-}(X)\,dK^g_{s+} \,\big|\,\cF_\tau\Big] \geq 0 ,
\end{align*}
which shows that $\barY_{\tau}\geq 0$ and thus $Y_{\tau}\geq \cE^g_{\tau,\sigma}(Y_\sigma)$. We have thus shown that $Y$
is a strong $\cE^g$-supermartingle.

\vskip 5 pt
\noindent {\it Step 2.} Our next goal is to show $Y$ can be characterized as the value process for a nonlinear optimal stopping problem associated with $\cE^g$ and $\zeta$.  Since $Y$ is a strong $\cE^g$-supermartingale and the nonlinear evaluation $\cE^g$ has the monotonicity property we have that $Y_\tau \geq \cE^g_{\tau,\sigma}(Y_\sigma) \geq \cE^g_{\tau,\sigma}(\zeta_\sigma)$ for all $\FF$-stopping times $\tau \leq \sigma$, which in turn implies that
\begin{align*}
Y_\tau \geq \sup_{\sigma}\mathcal{E}^g_{\tau,\sigma}(\zeta_\sigma).
\end{align*}
To show the reverse inequality, we fix $\tau \in \STF$ and we define the $\FF$-stopping time $\sigma^\varepsilon_t:= \inf\{s \geq
\tau : Y_s \leq \zeta_s + \varepsilon\}$. Since the obstacle $\zeta$ is upper-semicontinuous, using similar techniques to those in Section \ref{nsec3.4a}, one can show that $Y_{\sigma^\varepsilon_\tau} \leq \zeta_{\sigma^\varepsilon_\tau} + \varepsilon$ and $Y$ is an $\cE^g$-martingale on $\llb \tau,\sigma_\tau^\varepsilon \rrb$, that is, $Y$ is the solution to \eqref{yprime} on $\llb \tau,\sigma_\tau^\varepsilon \rrb$ with the terminal
condition $Y_{\sigma_\tau^\varepsilon}$. Using first the monotonicity property of $\cE^g$ stemming from Proposition \ref{npro1.1} and then
Proposition \ref{npro1.2x} with $g^1=g^2=g$, we deduce that for arbitrary $\varepsilon >0$
\begin{align*}
Y_{\tau}=\cE^g_{\tau, \sigma^\varepsilon_\tau}(Y_{\sigma^\varepsilon_\tau}) \leq \cE^g_{\tau, \sigma^\varepsilon_\tau}(\zeta_{\sigma^\varepsilon_t} + \varepsilon) \leq \cE^g_{\tau, \sigma^\varepsilon_\tau}(\zeta_{\sigma^\varepsilon_t} ) + C \varepsilon  \leq \sup_{\sigma\in \cT_{[\tau ,T]}(\FF)}\cE^g_{\tau, \sigma}(\zeta_{\sigma} ) + C \varepsilon .
\end{align*}
We have thus shown that for every $\tau, \sigma \in\STF$
\begin{align*}
Y_{\tau}=\sup_{\sigma\in \cT_{[\tau ,T]}(\FF)}\cE^g_{\tau, \sigma}(\zeta_{\sigma}),
\end{align*}
which means that $Y$ is the value process for a nonlinear optimal stopping problem associated with the nonlinear evaluation $\cE^g$
and the reward process $\zeta$.

\vskip 5 pt \noindent {\it Step 3.} Suppose that $h\leq g$ and the mapping $g$ is nondecreasing in $y$.  Then from the comparison property
of solutions to a GBSDE established in Proposition \ref{npro1.1}, we obtain for every $\tau \leq \sigma $,
\begin{align*}
\cE^{h}_{\tau, \sigma}(\wtzeta_{\sigma} ) \leq \cE^{g}_{\tau, \sigma}(\zeta_{\sigma} )
\end{align*}
and the assertion now follows by taking the supremum over all stopping times from $\STF$.
%\textcolor{red}{There is a step two in Grigorova's paper, see Theorem 5.3 on page 25 \cite{GIOQ2020}. We probably don't need the extension, but of course it can be included.}
\end{proof}

Let the mappings $f, \tilde{f} :\Omega \times [0,T]\times \RR \to \RR_+$ be given by $f(t,y) =(\eta_t-y)^+$ and $\tilde{f}(t,y) =-(y-\eta_t)^+$, respectively, where $\eta$ is a predetermined $\FF$-optional process. If we set $\gen (t,y)=nf(t,y)$ and $h(t,y)=(n+1)f(t,y)$ for $n\in\NN$ then $\gen$ and $h$ are nonnegative functions nonincreasing in $y$. Notice that $\tilde{\gen}(t,y)=n\tilde{f}(t,y)$ and $\tilde{h}(t,y)=(n+1)\tilde{f}(t,y)$
are nonpositive and nonincreasing in $y$. Then we have the following corollary to Proposition \ref{npro1.5}.

\bcor  \label{ncor1.3}
Assume that the mapping $f:\Omega \times [0,T]\times \RR \to \RR_+$ is nonnegative (resp., nonpositive) and such that $f(t,\cdot)$ is nonincreasing,
for every $t\in [0,T]$. For every $n\in\NN$, let the triplet $(Y^n,Z^n,K^n)$ be a unique solution to the reflected GBSDE
\begin{align} \label{neq1.12v}
Y^n_{\tau}=\zeta_T-\int_{\rrb\tau,T\rrb}Z^n_s\,dM_s+\int_{\rrb\tau,T\rrb}  nf(s,Y^n_s)\,d\newA_s+K^n_T-K^n_{\tau}
\end{align}
with the lower obstacle $\zeta$. Then the inequality $Y^{n+1}\geq Y^n$ (resp., $\wtY^{n+1}\leq \wtY^n$) holds for every $n\in\NN$.
\ecor

\begin{proof}
It suffices to apply Proposition \ref{npro1.5} to $\gen (t,y)= nf(t,y)$ and $h(t,y)= (n+1)f(t,y)$ where
$f$ is nonnegative (resp., nonpositive). It is obvious that $\gen$ is a nonincreasing function and $h(\omega,t,y)\ge \gen (\omega,t,y)$
(resp., $\gen (\omega,t,y)\geq h(\omega,t,y)$) since $g$ is a nonnegative (resp., nonpositive) function.
Hence Proposition \ref{npro1.5} implies that $Y^{n+1}\geq Y^n$ for every $n\in\NN$ (resp., $\wtY^{n+1}\leq \wtY^n$ where $(\wtY^n,\wtZ^n,\wtK^n)$
is a solution to \eqref{neq1.12v} with the generator $n\tilde{f}(s,Y^n_s)$).
\end{proof}

%%%%%%%%%%%%%%%%%%%%%%%%%%%%%%%%%%%%%%%%%%%%%%%%%%%%%%%%%%%%%%%%%%%%%%%%%%%%%%%%%%%%%%%%%%%%%%%%
\subsection{A Priori Estimates for Solutions to a Reflected GBSDE}    \label{nsec3.2}
%%%%%%%%%%%%%%%%%%%%%%%%%%%%%%%%%%%%%%%%%%%%%%%%%%%%%%%%%%%%%%%%%%%%%%%%%%%%%%%%%%%%%%%%%%%%%%%%

We first deal with {\it a priori} estimates for solutions to the reflected GBSDE.

\bp \label{npro1.1a}
Let $(Y^i,Z^i,K^i)\in \cS^2(\FF) \times \cH^2(M)\times \bcK$ for $i=1,2$ be a solution to the reflected GBSDE \eqref{RGBSDE} with data $(M,\newA,\gen^i,\zeta)$ for $i=1,2$, respectively, where the processes $A$ and $\zeta $ are bounded. Then for every $\beta>0$ there exists a constant $c>0$ such that the processes $Y=Y^1-Y^2$ and $Z=Z^1-Z^2$ satisfy
\begin{align}  \label{eqL6.11}
\EE\bigg[\sup_{t\in[0,T]}\ee_t|Y_t|^2+\int_{]0,T]} \ee_s|Z_s|^2\,d[M]_s\bigg] \leq  c\, \EE \bigg[\int_{]0,T]} \ee_s|\bar\gen_s|^2\,d\newA_s \bigg]
\end{align}
where $\ee_t:=e^{\beta \newA_t}$ and $\wh\gen_s:=\gen^1(s,Y^1_s)-\gen^2(s,Y^2_s)$.
\ep

\begin{proof}  The proof is similar to the proof of Proposition \ref{npro1.2}. Since $\beta >0$ we have that $d\ee_t:=\beta \ee_t\,d\newA_t$ and $1 \leq \ee_t \leq e^{\beta c_{\newA}}$ for every $t\in[0,T]$. Denote $K=K^1-K^2$ and recall that $K$ has a unique decomposition $K=K^r+K^g$ where $K^r$ is an $\FF$-adapted, c\`adl\`ag, nondecreasing process of finite variation and $K^g$ is an $\FF$-adapted, c\`agl\`ad, purely discontinuous, nondecreasing process of finite variation. Hence we obtain
\begin{align*}
Y_{\tau}=\int_{\rrb\tau,T\rrb} \bar\gen_s\,d\newA_s-\int_{\rrb\tau,T\rrb} Z_s\,dM_s+K^r_T-K^r_{\tau}+K^g_T-K^g_{\tau}.
\end{align*}
By applying the Gal'\v{c}uk-Lenglart formula (see, e.g., Theorem 8.2 in Gal'\v{c}uk \cite{G1981}) % or Lenglart (1980), Chapter VI, Section 3, page 538)
to $\ee_t|Y_t|^2$ and the Young inequality with a constant $\alpha>0$,  we obtain
\begin{align} \label{eqL6.12}
&\ee_t|Y_t|^2+\int_{]t,T]}\ee_s|Z_s|^2\,d[M]_s \nonumber\\
&=-\beta\int_{]t,T]}\ee_s|Y_s|^2\,d\newA_s+2\int_{]t,T]}\ee_sY_s\bar\gen_s\,dA_s-2\int_{]t,T]}\ee_sY_{s-}Z_s\,dM_s+2\int_{]t, T]}\ee_sY_{s-}\,dK^r_s \nonumber\\
&\ \ \ \  +  2\int_{[t, T[}\ee_sY_s\,dK^g_{s+} -\sum_{t<s\leq T}\ee_s(Y_s-Y_{s-})^2-\sum_{t\leq s< T}\ee_s(Y_{s+}-Y_{s})^2 \\
&\leq (\alpha^{-1}-\beta)\int_{]t,T]}\ee_s|Y_s|^2\,d\newA_s+\alpha\int_{]t,T]}\ee_s|\bar\gen_s|^2\,dA_s-2\int_{]t,T]}\ee_sY_{s-}Z_s\,dM_s +2\int_{]t, T]}\ee_sY_{s-}\,dK^r_s  \nonumber\\
&\ \ \ \   +2\int_{[t, T[}\ee_sY_s\,dK^g_{s+} -\sum_{t<s\leq T}\ee_s(Y_s-Y_{s-})^2-\sum_{t\leq s< T}\ee_s(Y_{s+}-Y_{s})^2. \nonumber
\end{align}
Next we show that $\int_{]t,T]}\ee_sY_{s-}\,dK^r_s \leq 0$ and $\int_{[t,T[}\ee_sY_{s}\,dK^g_{s+} \leq 0 $ for every $t\in [0,T]$.
We note that $K^r=(K^{1,c}-K^{2,c})+(K^{1,d}-K^{2,d})$. From the Skorokhod conditions in \eqref{RGBSDE} we obtain, for
all $s\in [0,T]$,
\begin{align*}
Y_{s}\,dK_{s}^{1, c}&=(Y_{s}^{1}-\zeta_{s})\,dK_{s}^{1,c}-(Y_{s}^{2}-\zeta_{s})\,dK_{s}^{1,c}\\
&\leq (Y_{s}^{1}- \bar\zeta_{s})\,dK_{s}^{1,c}-(Y_{s}^{2}-\zeta_{s})\,dK_{s}^{1, c} =-(Y_{s}^{2}-\zeta_{s})\,dK_{s}^{1,c}\leq 0
\end{align*}
where the last inequality holds since the process $K^{c}$ is nondecreasing. By symmetry, we obtain $Y_{s}\, dK_{s}^{2, c} \geq 0$.
Furthermore, for any $s\in [0,T]$,
\begin{align*}
Y_{s-}\Delta K_{s}^{1,d}&=(Y_{s-}^{1}-\zeta_{s})\Delta K_{s}^{1,d}-(Y_{s-}^{2}-\zeta_{s}) \Delta K_{s}^{1,d}\\
&\leq (Y_{s-}^{1}-\bar\zeta_{s})\Delta K_{s}^{1,d}-(Y_{s-}^{2}-\zeta_{s}) \Delta K_{s}^{1,d} =-(Y_{s-}^{2}-\zeta_{s}) \Delta K_{s}^{1,d} \leq 0
\end{align*}
and $Y_{s-} \Delta K_{s}^{2,d} \geq 0$. Similarly, $\int_{[t, T[}\ee_s Y_{s}\,dK^g_{s+}=\sum_{t\leq s <T}\ee_s Y_s\Delta^+K^g_s$ for all $s\in [0,T]$.
Note that
\begin{align*}
Y_s\Delta^+K^g_s=(Y_s^1-Y_s^2)\Delta^+K^{1,g}_s-(Y_s^1-Y_s^2)\Delta^+K^{2,g}_s
\end{align*}
and, for all $s\in [0,T]$,
\begin{align*}
(Y_s^1-Y_s^2)\Delta^+K_{s}^{1,g}&=(Y_{s}^{1}-\zeta_{s})\Delta^+K_{s}^{1,g}-(Y_{s}^{2}-\zeta_{s})\Delta^+K_{s}^{1,g}
=-(Y_{s}^{2}-\zeta_{s})\Delta^+K_{s}^{1,g} \leq 0
\end{align*}
and $(Y_s^1-Y_s^2)\Delta^+K^{2,g}_s \geq 0$. Hence \eqref{eqL6.12} gives, for any $\beta>\alpha^{-1}$ and $t\in [0,T]$,
\begin{align}\label{eqL6.13}
\ee_t|Y_t|^2+\int_{]t,T]}\ee_s|Z_s|^2\,d[M]_s \leq \alpha\,\int_{]t,T]}\ee_s|\bar\gen_s|^2\,d\newA_s-2\int_{]t,T]}\ee_sY_{s-}Z_s\,dM_s
\end{align}
and thus, by taking the expectation on both sides, we obtain
\begin{align}\label{eqL6.14}
\EE\bigg[\int_{]t,T]}\ee_s|Z_s|^2\,d[M]_s\bigg] \leq \alpha\,\EE\bigg[\int_{]t,T]}\ee_s|\bar\gen_s|^2\,d\newA_s\bigg].
\end{align}
In addition, taking the essential supremum and expectation in \eqref{eqL6.13} gives
\begin{align} \label{eqL6.15}
\EE\bigg[\esssup_{\tau \in \STF} \ee_\tau|Y_\tau|^2\bigg]\leq \EE \bigg[\alpha\,\int_{]0,T]} \ee_s|\bar\gen_s|^2\,d\newA_s\bigg]
 +2\,\EE\bigg[\esssup_{\tau \in \STF}\Big|\int_{]0,\tau]}\ee_\tau Y_{\tau-}Z_\tau\,dM_\tau\Big|\bigg].
\end{align}
An application of the Burkholder--Davis--Gundy inequality with $p=1$ similar to \eqref{eq1.16d} yields
\begin{align} \label{eqL6.16}
2\,\EE\bigg[\esssup_{\tau \in \STF} \Big|\int_{]0,\tau]}\ee_sY_{s-}Z_s\,dM_s\Big|\bigg]
\leq \frac{1}{4}\,\EE \bigg[\esssup_{\tau \in \STF}\ee_\tau|Y_\tau|^2\bigg]
+4c_1^2\,\EE\bigg[\int_{]0,T]} \ee_s|Z_s|^2\,d[M]_s\bigg]
\end{align}
where the constant $c_1$ is independent of $\alpha,\beta$ and the last inequality holds since $2ab\leq a^2+b^2$ for all real numbers $a,b$.
By combining \eqref{eqL6.14}, \eqref{eqL6.15} and \eqref{eqL6.16}, we obtain
\begin{align*}
\frac{3}{4}\,\EE\bigg[\esssup_{\tau\in\STF}\ee_\tau|Y_\tau|^2\bigg]\leq (1+4c_1^2)\alpha\,\EE\bigg[\int_{]0,T]} \ee_s|\bar\gen_s|^2\,d\newA_s \bigg]
\end{align*}
and, finally,
\begin{align*}
\EE\bigg[\esssup_{\tau\in\STF}\ee_\tau|Y_\tau|^2+\int_{]0,T]}\ee_s|Z_s|^2\,d[M]_s\bigg]\leq c\,\EE\bigg[\int_{]0,T]} \ee_s|\bar\gen_s|^2\,d\newA_s \bigg]
\end{align*}
where $c=\frac{7+16c_1^2}{3}\alpha$. 
\end{proof}

Assume that the generator $g^1$ is Lipschitz continuous with a constant $L$. Then
\begin{align*}
|\bar\gen_s|^2 &= |\gen^1(s,Y^1_s)-\gen^1(s,Y^2_s) +\gen^1(s,Y^2_s) -\gen^2(s,Y^2_s)|^2 \\
&\leq  2|\gen^1(s,Y^1_s)-\gen^1(s,Y^2_s)|^2 + 2|\gen^1(s,Y^2_s) -\gen^2(s,Y^2_s)|^2 \\
&\leq  2L^2 |Y^1_s-Y^2_s|^2 + 2 |\wh\gen_s|^2
\end{align*}
where we denote $\wh\gen_s = \gen^1(s,Y^2_s) -\gen^2(s,Y^2_s)$. Consequently, in \eqref{eqL6.12} we get
\begin{align*}
 (\alpha^{-1}-\beta)\int_{]t,T]}\ee_s|Y_s|^2\,d\newA_s+\alpha\int_{]t,T]}\ee_s|\bar\gen_s|^2\,dA_s \leq
 c_{\alpha,\beta,L} \int_{]t,T]}\ee_s|Y_s|^2\,d\newA_s+ 2 \alpha\int_{]t,T]}\ee_s |\wh\gen_s|^2 \,dA_s .
\end{align*}
where $c_{\alpha,\beta,L} =  (\alpha^{-1}+ 2 \alpha L^2 -\beta)$. Then, by arguing as in the proof of Proposition \ref{npro1.1a}, we obtain  the following
modification of \eqref{eqL6.13}, which is valid for every $\beta > \alpha^{-1}+ 2 \alpha L^2$
\begin{align} \label{eqL6.13c}
\ee_t|Y_t|^2+\int_{]t,T]}\ee_s|Z_s|^2\,d[M]_s \leq \alpha\,\int_{]t,T]}\ee_s|\wh \gen_s|^2\,d\newA_s-2\int_{]t,T]}\ee_sY_{s-}Z_s\,dM_s .
\end{align}
By taking the conditional expectation of both sides with respect to $\cF_t$ we get
\begin{align*}
\ee_t|Y_t|^2+ \EE\bigg[\int_{]0,T]}\ee_s|Z_s|^2\,d[M]_s\,\Big|\, \cF_t \bigg]\leq \EE\bigg[\int_{]t,T]} \ee_s|\wh \gen_s|^2\,d\newA_s\,\Big|\, \cF_t \bigg],
\end{align*}
which holds for a sufficiently large $\beta$. If, in addition, an increasing process $A$ is bounded, we obtain the existence
of a constant $c$ such that for all $t\in [0,T]$
\begin{align*}
|Y_t|^2 \leq c\, \EE\bigg[\sup_{s\in[t,T]}|\wh \gen_s|^2\,\Big|\, \cF_t\bigg].
\end{align*}
In particular, if $|\gen^1(s,y)-\gen^2(s,y)|\leq \varepsilon $ for all $(s,y)$, then $\sup_{t\in[0,T]}|Y^1_t-Y^2_t|\leq c \varepsilon$ for every $t\in [0,T]$.

%%%%%%%%%%%%%%%%%%%%%%%%%%%%%%%%%%%%%%%%%%%%%%%%%%%%%%%%%%%%%%%%%%%%%%%%%%%%%%%%%%%%%%%%%%%%%%%%
\subsection{Existence of a Solution to a Reflected GBSDE}         \label{nsec3.3}
%%%%%%%%%%%%%%%%%%%%%%%%%%%%%%%%%%%%%%%%%%%%%%%%%%%%%%%%%%%%%%%%%%%%%%%%%%%%%%%%%%%%%%%%%%%%%%%%

We now examine the existence of a solution $(Y,Z,K)$ to the reflected GBSDE \eqref{RGBSDE}. For brevity, we only work under the assumption that
the driver $A$ is bounded but we conjecture that the proof can be extended to the case of a non-negative, nondecreasing generator and square-integrable driver by proceeding analogously to the proof of Proposition \ref{npro1.4}.
We first establish the existence and uniqueness of a solution to the reflected GBSDE \eqref{RGBSDE} when the generator $g$ does not depend on $Y$.

\bl \label{nlem1.2}
Let Assumptions \ref{nass1.1} and \ref{nass1.2} hold and $g(s,\omega,y)=g(s,\omega):=g_s$ be a fixed process. 
If the processes $A$ and $\zeta $ are bounded, then the reflected GBSDE \eqref{RGBSDE} with data $(M,\newA,\gen,\zeta)$  admits a unique solution $(Y,Z,K)\in \cS^2(\FF) \times \cH^2(M)\times \bcK$.
\el

\begin{proof}
For $\nu \in \STF$, we first denote
\begin{align} \label{eqL6.17}
Y(\nu) := \esssup_{\tau \in \STFnu}\EE \bigg[ \zeta_\tau+\int_{\rrb \nu, \tau\rrb} g_s\,d\newA_s\,\Big |\,\mathcal{F}_\nu\bigg].
\end{align}
By the classical aggregation result (see, e.g., Theorem 8.2 in Grigorova et al. \cite{GIOQ2020}) there exists a l\`adl\`ag, right-upper-semicontinuous
process $Y$ such that, for all $\nu \in \STF$,
\begin{align} \label{eqL6.18}
Y_\nu=Y(\nu)
\end{align}
and the process $(Y_t+\int_{]0,t]} g_s\,d\newA_s)_{t \in [0,T]}$ is the smallest strong supermartingale that dominates the process
$[\zeta_.+\int_{]0,\cdot]} g_s\, d\newA_s]$.

\noindent {\it Step 1.} We show that $Y\in S^2{(\FF)}$ is a candidate for the first component in the solution $(Y,Z,K)$ to the reflected GBSDE \eqref{RGBSDE}.
An application of Jensen's inequality gives, for all $\nu \in \STF$,
\begin{align*}  % \label{eqL6.19}
\left|Y_{\nu}\right| & \leq \underset{\tau \in  \STFnu}{\operatorname{ess} \sup }
\EE \bigg[\big|\zeta_{\tau}\big|+\big|\int_{\rrb \nu, \tau\rrb}  g_s \,d\newA_s\big| \,\Big|\,\mathcal{F}_{\nu}\bigg] \\
& \leq \EE\bigg[\underset{\tau \in  \STF}{\operatorname{ess} \sup }\,|\zeta_{\tau}|+\int_{]0,T]}|g_s|\, d\newA_s \,\Big|\,\mathcal{F}_{\nu}\bigg]=
\EE [X|\,\mathcal{F}_{\nu}]
\end{align*}
where $X:=\esssup_{\tau \in \STF}|\zeta_{\tau}|+\int_{]0,T]}|g_s|\,d\newA_s$. The Cauchy-Schwarz inequality yields
\begin{align*}
\EE[X^2]\leq c\,||\zeta||^2_{S^2(\FF)}+cT\,\EE \bigg[\int_{]0,T]} |g_s|^2\,d\newA_s\bigg]
\end{align*}
where $c>0$ is a positive constant, which may vary from line to line.
Since $\esssup_{\nu \in \STF}|Y_\nu|^2\leq \sup_{t\in [0,T]}|\EE[X\,|\,\mathcal{F}_t]|^2$ the Doob martingale inequality in $L^2$ gives
\begin{align*}
\EE\Big[\esssup_{\nu \in \STF}|Y_\nu|^2\Big]\leq \EE \Big[\sup_{t\in [0,T]}|\EE[X\,|\,\mathcal{F}_t]|^2\Big]\leq c \,\EE[X^2]\leq  c\,||\zeta||^2_{S^2(\FF)}
+cT\,\EE \bigg[\int_{]0,T]} |g_s|^2\,d\newA_s \bigg]
\end{align*}
where the penultimate inequality is due to the right-continuity of the process $\EE[X|\mathcal{F}_t]_{t \in [0,T]}$.

\noindent {\it Step 2.} We will show that there exists a pair $(Z,K) \in (\cH^2(M),\bcK)$ such that $(Y,Z,K)$ satisfies the reflected GBSDE and $\zeta$ is indeed an obstacle. Since $(Y_t+\int_{]0,t]} g_s\,d\newA_s)_{t \in [0,T]}$ is a strong supermartingle and $g$ satisfies Assumption \ref{nass1.2} by the Doob-Meyer-Mertens decomposition (see Theorem \ref{nth1.1}) we obtain
\begin{align}\label{eqL6.20}
Y_t=-\int_{]0,t]} g_s\,d\newA_s+N_t-K^c-K^d-K^g
\end{align}
where $N=Z\bigcdot M$ is a uniformly integrable martingale by the predictable representation property of $M$ and $K=K^c+K^d+K^g$ is a l\`adl\`ag, nondecreasing, $\FF$-predictable process satisfying the Skorokhod conditions in \eqref{RGBSDE} (see \cite{EK1981}). Recall that $Y_T=Y(T)=\zeta_T$. From
\eqref{eqL6.17} and \eqref{eqL6.18}, we obtain $Y_\nu=Y(\nu) \geq \zeta_{\nu}$ for every $\nu \in \STF$. Therefore, by a classical result of Theorem IV.84 in \cite{DM1975}, we conclude that $Y_t\geq \zeta_t, \, 0\leq t \leq T$, a.s.

\noindent {\it Step 3.} Now we establish the uniqueness of a solution to the reflected GBSDE \eqref{RGBSDE}. Let $(Y',Z')$ be another solution to
\eqref{RGBSDE}. By Proposition \ref{npro1.1a} with $g^1(\cdot)=g^2(\cdot)=g$, we obtain $Y=Y'$ in $\cS^2(\FF)$ and $Z=Z'$ in $ \cH^2(M)$. Then the uniqueness of $K$ follows from the uniqueness of the Doob-Meyer-Mertens decomposition.
\end{proof}

We are in a position to prove the existence and uniqueness of the solution to reflected GBSDE \eqref{RGBSDE} with a general generator $g$.

\bp \label{npro1.6}
Let Assumptions \ref{nass1.1} and \ref{nass1.2} hold with a bounded process $A$ and $\zeta$ be an $\FF$-optional, nonnegative, bounded process.
Then the reflected GBSDE \eqref{RGBSDE} with data $(M,\newA,\gen,\zeta)$  admits a unique solution $(Y,Z,K)\in \cS^2(\FF) \times \cH^2(M)\times \bcK$.
\ep

\begin{proof}
We will extend the proof of Proposition \ref{npro1.3}. Note that $Y\in S^2(\FF)$ in \eqref{RGBSDE} is not necessarily right-continuous and we equipped $\cS^2$ with the norm $||Y||^2_{\cS^2_\beta(\FF)}=\EE[\esssup_{\tau\in \STF}\ee_\tau |Y_\tau|^2]$  under which $\cS^2$ is still a Banach space. We observe that for $\beta > 0$ the norms $||\cdot||_{\cS^2(\FF)}$ and $||\cdot||_{\cS^2_\beta(\FF)}$ are equivalent on $\cS^2(\FF)$ and denote $\cS_{\beta}^2(\FF)$ the space $\cS^2(\FF)$ endowed with the norm $||\cdot||_{\cS^2_\beta(\FF)}$. Next we adopt the Banach fixed point theorem.

Let the mapping $\Phi : \cS_{\beta}^2(\FF) \to \cS_{\beta}^2(\FF)$ be defined as follows: for any given $\ww \in \cS_{\beta}^2(\FF)$ we set $\Phi(\ww):=Y^{\ww}$ where the triplet $(Y^{\ww},Z^{\ww},K^{\ww})$ is a unique solution to the reflected GBSDE (see Lemma \ref{nlem1.2})
\begin{align} \label{eqL6.21}
Y^{\ww}_t=\zeta_T-\int_{]t,T]}Z^{\ww}_s\,dM_s+\int_{]t,T]} g(s,\ww_s)\,d\newA_s+K^{\ww}_T-K^\ww_t
\end{align}
where $Y^{\ww}\geq \zeta$ and the fixed generator is independent of $Z^{\ww}$.

Our next goal is to demonstrate that there exists a unique process $\wh\ww \in \cS_{\beta}^2(\FF) $ such that $\Phi (\wh\ww)=\wh\ww$.
Then the corresponding process $\wh z \in \cH^2(M)$ and $\wh k\in \bcK$ can be determined by \eqref{eqL6.21}, that is, we would obtain
\begin{align*}
\wh\ww_t=\zeta_T-\int_{]t,T]}\wh z_s\,dM_s+\int_{]t,T]} g(s,\wh\ww_t)\,d\newA_s+\wh k_T-\wh k_t
\end{align*}
with $\wh k_0=0$ and the nondecreasing process $\wh k$ satisfies the Skorokhod conditions in \eqref{RGBSDE}.
It is clear that it suffices to show that the mapping $\Phi:\cS_{\beta}^2(\FF)\to\cS_{\beta}^2(\FF)$ is a contraction for a sufficiently large $\beta$.
We take $\ww', \ww''\in \cS_{\beta}^2(\FF)$ and denote $Y^{\ww'}=\Phi(\ww')$ and $Y^{\ww''}=\Phi(\ww'')$. For the simplicity of notation, we write $\yyy:=Y^{\ww'}-Y^{\ww''}=\Phi(\ww')-\Phi(\ww''),\, \zzz :=Z^{\ww'}-Z^{\ww''},\,k:=K'-K''$ and $\ww:= \ww'- \ww''$.
It is clear from \eqref{neq1.17} that $\yyy$ satisfies the reflected GBSDE
\begin{align*}
\yyy_t=-\int_{]t,T]}\zzz_s\,dM_s+\int_{]t,T]} \big(g(s,\ww'_s)-g(s,\ww''_s)\big)\,d\newA_s+k_T-k_t
\end{align*}
where $|g(s,\ww'_s)-g(s,\ww''_s)|\leq L|\ww_s|$ since $g(s,\cdot)$ is Lipschitz continuous with a constant $L$. By applying
Proposition \ref{npro1.1a}, we obtain
\begin{align*}
\EE\bigg[\esssup_{\tau \in \STF} \ee_\tau |\ww_\tau |^2 +\int_{]0,T]} \ee_s|z_s|^2\,d[M]_s\bigg] \leq \alpha L^2\frac{7+16c_1^2}{3}
\,\EE \bigg[\int_{]0,T]} \ee_s|w_s|^2\,d\newA_s \bigg]
\end{align*}
where $\alpha>0$, $\beta>\frac{1}{\alpha}$ and a constant $c_1>0$ is independent of $\alpha, \beta, L$ and thus
\begin{align*}
\EE\bigg[\esssup_{\tau \in \STF} \ee_\tau |\ww_\tau |^2 +\int_{]0,T]} \ee_s|z_s|^2\,d[M]_s\bigg] \leq \alpha c_A L^2 (T+1)\frac{7+16c_1^2}{3}
\, \EE \bigg[ \esssup_{\tau \in \STF} \ee_\tau |\ww_\tau |^2 \bigg]
\end{align*}
where $c_A>0$ is such that $\newA_T \leq c_A$. Thus we conclude that $\Phi$ is a contraction when $0< \alpha<c_A^{-1} L^{-2}\frac{3}{(7+16c_1^2)(T+1)}$ and $\beta>\alpha^{-1}$ and thus there exists a unique solution $(Y,Z,K) \in  \cS^2(\FF) \times \cH^2(M)\times \bcK$ to the reflected GBSDE \eqref{RGBSDE}.
\end{proof}

%%%%%%%%%%%%%%%%%%%%%%%%%%%%%%%%%%%%%%%%%%%%%%%%%%%%%%%%%%%%%%%%%%%%%%%%%%%%%%%%%%%%%%%%%%%%%%%%
\subsection{Penalization Schemes for a Reflected GBSDE}             \label{nsec3.4a}
%%%%%%%%%%%%%%%%%%%%%%%%%%%%%%%%%%%%%%%%%%%%%%%%%%%%%%%%%%%%%%%%%%%%%%%%%%%%%%%%%%%%%%%%%%%%%%%%

To formulate results on penalization schemes for the reflected GBSDE, we recall that $\barS = S^r\cup\{T\}$ where $S^r=S^r(A)$ is the right support of the process $A$. The first  penalization scheme is analogous to the case of a penalized GBSDE, which was examined in Section \ref{nsec3.4}. Similar to
Theorem \ref{th3.1}, we will show in Theorem \ref{th6.1} that the process $Y=\lim_{\,n\to \infty} Y^n$ can be interpreted as the value process of an optimal stopping problem with the reward process $\gamma $ given by
\begin{align}  \label{eq8n4.2}
\gamma_t:=\zeta_T\I_{\{t=T\}}+(\zeta\vee \eta \I_{\barS})_t\I_{\{t<T\}}.
\end{align}
We henceforth postulate that Assumptions \ref{nass1.1} and \ref{nass1.2} are satisfied with a bounded process $A$.

\bt \label{th6.1}
Let the $\FF$-optional and bounded process $\zeta$ and (resp., the $\FF$-optional and bounded process $\eta$) be l\`adl\`ag and right-upper-semicontinuous (resp., right-continuous). Consider the sequence of solutions $(Y^n,Z^n,K^n)$ to the reflected GBSDE
\begin{align} \label{eq6.6}
Y^n_{\tau}=\zeta_T-\int_{\rrb\tau,T\rrb} Z^n_s\,dM_s+\int_{\rrb\tau,T\rrb} n(\eta_s-Y^n_s)^+\,dA_s+K^n_T-K^n_{\tau}
\end{align}
where an $\FF$-adapted, l\`adl\`ag, nondecreasing process $K^n$ satisfies the Skorokhod conditions with the lower obstacle $\zeta$.
Then the sequence $Y^n$ converges monotonically to a process $Y$, which satisfies, for every $t\in [0,T]$,
\begin{equation} \label{eq6.5}
Y_t=\esssup_{\sigma\in\STFt}\EE\big[\gamma_\sigma|\,\cF_t\big]
\end{equation}
where the process $\gamma$ is given by \eqref{eq8n4.2}. In addition, the triplet $(Y,Z,K)=\lim_{\,n \to \infty}(Y^n,Z^n,K^n)$ is a unique solution to the reflected BSDE
\begin{align}  \label{eq8.4.2}
Y_t=\zeta_T-\int_{]t ,T]}Z_s\,dM_s+K_T-K_{t}
\end{align}
where an $\FF$-adapted, l\`adl\`ag, nondecreasing process $K$ satisfies the Skorokhod conditions with the lower obstacle $\zeta $.
\et

\begin{proof}
We start by noticing that the existence of a unique solution $(Y^n,Z^n,K^n) \in \cS^2(\FF) \times \cH^2(M)\times \bcK$ to the reflected GBSDE
\eqref{eq6.6} follows from Proposition \ref{npro1.6}. Furthermore, the sequence $Y^n$ of processes is monotonically increasing as $n\rightarrow \infty$ (see Corollary \ref{ncor1.3}) and the limit $Y=\lim_{\,n\rightarrow \infty}Y^n$ is well defined.

\noindent {\it Step 1.}  Our first goal is to show that, for every $n\in\NN$,
\begin{align}  \label{eq6.6a}
Y^n_t=\esssup_{\sigma\in\STFt}\EE\big[Y^n_\sigma \wedge \gamma_\sigma\,|\,\cF_t\big].
\end{align}
% where $\gamma_t := (\zeta\vee \eta\I_{\barS})_t\I_{\{t<T\}}+\zeta_T \I_{\{t=T\}}$.
To prove \eqref{eq6.6a}, we fix $n$ and we observe that $Y^n$ is a supermartingale and thus, for every $\sigma\in\STFt$,
\begin{align} \label{eq6.6b}
Y^n_t \geq\EE\big[Y^n_\sigma\,|\,\cF_t\big]\geq \EE\big[Y^n_\sigma \wedge \gamma_\sigma\,|\,\cF_t\big]
\end{align}
where the second inequality is obvious.
To show the reverse inequality, we fix $t\in [0,T[$ and we define $\nu =\sigma_t^n\wedge\tau_t^n \in \STFt$ where for any fixed $\delta>0$ we set  (as usual, $\inf \emptyset =T$)
\begin{align*}
\sigma_t^n:=\inf\{s\in [t,T]: Y^n_s \leq \zeta_s + \varepsilon \}, \quad \tau_t^n:=\inf\{s\in [t,T]:\wtK^n_s-\wtK^n_t>0\},\
\end{align*}
where $\varepsilon :=0.5(Y^n_t-\zeta_t)\delta $ and the continuous, nondecreasing process $\wtK^n$ is given by $\wtK^n_t:=\int_{]0,t]} n(\eta_s-Y^n_s)^+\,dA_s$. We will check that $Y^n_\nu = Y^n_\nu\wedge ( (\zeta+\varepsilon)\vee \eta\I_S)_\nu$ on the event $\{\nu <T\} =\{\sigma_t^n\leq \tau_t^n<T\} \cup \{\tau_t^n< \sigma_t^n\}=E_1 \cup E_2$. It is obvious that $Y^n_\nu = \zeta_T$ on the event $E_3:=\{\nu =T\}$.

On the event $E_1=\{t<\sigma_t^n\leq \tau_t^n<T\}$, we have $Y^n_{\sigma_t^n-}-\zeta_{\sigma_t^n-} \geq \varepsilon$ and thus $\Delta K^{n,d}_{\sigma_t^n}=0$, which implies that $Y^n$ is a martingale on $\llb t,\sigma_t^n\rrb$. Furthermore, if $\Delta^+ K^{n,g}_{\sigma_t^n}>0$, then the Skorokhod condition
gives $Y^n_{\sigma_t^n}=\zeta_{\sigma_t^n}$ and if $\Delta^+ K^{n,g}_{\sigma_t^n}=0$, then $Y^n$ is continuous at $\sigma_t^n$ and  $\zeta_{\sigma_t^n}\leq
Y^n_{\sigma_t^n}\leq \zeta_{\sigma_t^n}+\varepsilon$ since $\zeta$ is assumed to be right-upper-semicontinuous. We conclude that on $E_1$
we have $Y^n_\nu = Y^n_\nu  \wedge (\zeta_\nu + \varepsilon) = Y^n_\nu\wedge ( (\zeta+\varepsilon) \vee \eta\I_S)_\nu$ where the second equality is a trivial consequence of the first one.

On the event $E_2=\{\tau_t^n < \sigma_t^n\}$, the process $Y^n$ is right-continuous at $\tau_t^n$ and hence from the definition of $\tau_t^n$ we obtain $Y^n_{\tau_t^n}= Y^n_{\tau_t^n+} \leq \eta_{\tau_t^n+} = \eta_{\tau_t^n}$ where the last inequality follows from the right-continuity of $\eta$.
We note also that the $\FF$-stopping time $\tau_t^n$ has values in $S$ so that $\eta_{\tau_t^n}=(\eta\I_S)_{\tau_t^n}$ and thus we have $Y^n_\nu = Y^n_\nu  \wedge (\eta\I_S)_\nu = Y^n_\nu\wedge ((\zeta+\varepsilon)\vee \eta\I_S)_\nu$ on $E_2$ where the second equality is obvious.  It is also clear that $Y^n$ is a martingale on $\llb t,\tau_t^n\rrb$ since the continuous, nondecreasing process $\wt K^n$ and the nondecreasing process $K^n$ are constant on that interval.

Recall that $\varepsilon=0.5(Y^n_t-\zeta_t)\delta $ and the processes $Y^n$ and $\zeta$ are bounded so that $\varepsilon \leq c \delta$ for some constant $c$. Let us denote $\gamma^{\varepsilon}_t:=((\zeta+\varepsilon)\vee \eta\I_S)_t\I_{\{t<T\}}+\zeta_T\I_{\{t=T\}}$. Since $Y^n$ is a martingale on $\llb t,\nu \rrb$ we have
\begin{align} \label{eq6.6c}
Y^n_t=\EE\big[Y^n_\nu\,|\,\cF_t\big]= \EE\big[Y^n_\nu \wedge \gamma^{\varepsilon}_\nu\,|\,\cF_t\big]
 \leq \EE\big[Y^n_\nu \wedge \gamma_\nu\,|\,\cF_t\big] + c\delta \leq \EE\big[Y^n_\nu \wedge \gamma_\nu\,|\,\cF_t\big]
\end{align}
where the last inequality holds since $\delta$ is any positive number. By combining \eqref{eq6.6b} with \eqref{eq6.6c} we conclude that
\eqref{eq6.6a} is satisfied for every $n\in\NN$.

\noindent {\it Step 2.} We are now ready to show that  \eqref{eq6.5} is valid. For any $\tau \in \bSTFt$ equation \eqref{eq6.6} gives
\begin{align*}
Y^n_{\tau} & = \zeta_{T}-\int_{\rrb\tau,T\rrb} Z^n_s\,dM_s+\int_{\rrb \tau,T\rrb} n(\eta_s-Y^n_s)^+\,dA_s+K_{T}-K_{\tau}
\end{align*}
and, by applying the comparison theorem for a generalized GBSDE (see Proposition \ref{npro1.5x}),
we obtain the inequality $Y^n \geq \whY^n $ where $(\whY^n,\whZ^n,\whK^n)$ solves the linear reflected GBSDE
\begin{align*}
\whY_\tau^n&=\zeta_T-\int_{\rrb\tau,T\rrb}\wh  Z^n_s\,dM_s+\int_{\rrb \tau,T\rrb}n(\eta_s-\whY^n_s)\,dA_s+\whK^n_T-\whK^n_\tau\\
&=\zeta_T-\int_{\rrb\tau,T\rrb} \whZ^n_s\,dM_s+ \int_{\rrb \tau, T\rrb} n(\eta_s- \whY^n_s)\,dA_s+\whK^n_T-\whK^n_\tau\\
&\geq  \zeta_T-\int_{\rrb\tau,T\rrb} \barZ^n_s\,dM_s+ \int_{\rrb \tau, T\rrb} n(\eta_s- \barY^n_s)\,dA_s = \barY^n_{\tau}
\end{align*}
where the inequality holds since the generator $g^n(t,y)=n(\eta_t-y)$ is linear and the process $\whK$ is nondecreasing on $\llb \tau, T\rrb$.
Furthermore, by solving the linear GBSDE
\begin{align*}
\barY^n_{\tau} = \zeta_T-\int_{\rrb \tau,T\rrb} \barZ^n_s\,dM_s+ \int_{\rrb \tau,T\rrb} n(\eta_s-\barY^n_s)\,dA_s 
\end{align*}
 we obtain $Y^n_{\tau} \geq \whY^n_{\tau} \geq \barY^n_{\tau}$ where
\begin{align*}
\barY_{\tau}^n=\EE\big[\zeta_{T}\cE_{\tau,T}(-A^n)+(\I_{\rrb \tau,T\rrb}\eta\cE_{\tau,\cdot}(-A^n)\bcd A^n)_{T}\, | \cF_\tau\big]	
\end{align*}
where the sequence of random variables $\cE_{\tau,T}(-A^n)$ converges to $\I_{\{\tau=T\}}$ as $n\rightarrow \infty$ and thus, by Lemma \ref{lem3.4x} and the right-continuity of $\eta$, we obtain
\begin{align} \label{VR1}
Y_\tau=\lim_{\,n\rightarrow \infty}Y^n_{\tau} \geq \lim_{\,n \rightarrow \infty} \barY^n_{\tau}=\zeta_{\tau}\I_{\{\tau= \sigma\}}+\eta_{\tau}\I_{\{\tau< \sigma\}}.
\end{align}
Using the fact that $Y\geq \barY^n \geq 0$ and $Y \geq \zeta$, we deduce that for any stopping time $\sigma\in\STFt$, on the event $\{\sigma < T\}$, we have from \eqref{VR1}
\begin{align*}
Y_\sigma \geq \zeta_{\sigma} \vee (\eta\I_{\barS})_\sigma
\end{align*}
and $Y_T = \zeta_T$ on the event $\{\sigma = T\}$. Consequently,
\begin{align*}
Y_\sigma \geq \zeta_{T}\I_{\{\sigma = T\}}+\zeta_{\sigma} \vee (\eta\I_{\barS})_\sigma  \I_{\{\sigma < T\}} = \gamma_{\sigma}.
\end{align*}
Since we clearly have $Y^n_t \leq \esssup_{\sigma\in\STFt}\EE\big[\gamma_\sigma\,|\,\cF_t\big]$, it suffices to show that
\begin{align*}
 Y_t  \geq  \esssup_{\sigma \in\STFt}\EE\big[\gamma_\sigma\,|\,\cF_t\big].
\end{align*}
The above inequality follows from the observation that, by the monotone convergence theorem, the process $Y$ is a supermartingale dominating the reward process $\gamma$ and the minimality property of the Snell envelope. We thus conclude that \eqref{eq6.5} holds.

Finally, the representation of $Y$ as the solution of the reflected BSDE \eqref{eq8.4.2} follows by noticing that the process $\gamma$ is again upper-semicontinuous and applying the classic Doob-Meyer-Mertens decomposition to \eqref{eq6.5}.
\end{proof}

We conclude this work by examining another penalization scheme for a reflected GBSDE where the process $\zeta$ is still the lower obstacle but the role of the process $\eta$ differs from the previous result and thus the limiting process $\wtY$ is expected to represent the value of a Dynkin game and hence also to be a solution to a doubly reflected GBSDE (see Remark \ref{remx.x}).

\bt  \label{th6.2}
Let the $\FF$-optional and bounded process $\zeta$ (resp., the $\FF$-optional and bounded process $\eta$) be right-upper-semicontinuous (resp., right-continuous). Consider the sequence of solutions $(\wtY^n,\wtZ^n,\wtK^n)$ to the reflected GBSDE
\begin{align*}
\wtY^n_{\tau}=\zeta_T-\int_{\rrb\tau,T\rrb}\wtZ^n_s\,dM_s-\int_{\rrb \tau, T\rrb} n(\wtY^n_s-\eta_s)^+\,dA_s+\wtK^n_T-\wtK^n_{\tau}
\end{align*}
where  $\wtK^n$ satisfies the Skorokhod conditions with the lower obstacle $\zeta$. If the inequality $\zeta \leq \eta$ holds on the set $\barS$, then the sequence $\wtY^n$ converges monotonically to the process $\wtY$ which satisfies, for every $t\in [0,T]$,
\begin{equation} \label{eq6.7}
\wtY_t=\essinf_{\tau\in\bSTFt}\esssup_{\sigma\in\STFt}\EE(\Theta(\sigma,\tau)\,|\,\cF_t)=
\esssup_{\sigma\in\STFt}\essinf_{\tau\in\bSTFt}\EE(\Theta(\sigma,\tau)\,|\,\cF_t)
\end{equation}
where $\Theta(\sigma ,\tau ):=\zeta_{\sigma}\I_\seq{\tau > \sigma}+(\zeta\vee \eta)_\tau \I_\seq{\tau \leq \sigma}$.
\et

\begin{proof}
Let $(\wtY^n,\wtZ^n,\wt{K}^n)$ be the unique solution  in $\cS^2(\FF) \times \cH^2(M)\times \bcK$ to the reflected GBSDE (see Proposition \ref{npro1.6})
\begin{align}  \label{eq6.9m}
\wtY^n_{\tau}=\zeta_T-\int_{\rrb\tau,T\rrb} \wtZ^n_s\,dM_s-\int_{\rrb\tau,T\rrb} n(\wtY^n_s-\eta_s)^+\,dA_s+\wtK^n_T-\wtK^n_{\tau}
\end{align}
where $\wtY^n \geq \zeta$ and the Skorokhod conditions are satisfied by an $\FF$-adapted, nondecreasing process $\wtK^n$.
%Observe that, as opposed to the maximizer's case, the recovery process $\eta$ is acting here similarly to an upper obstacle for $\wtY^n$.
We note that the sequence $\wtY^n$ is monotonically decreasing as $n\rightarrow \infty$ (see Corollary \ref{ncor1.3}) and the limit $\wtY=\lim_{\,n\rightarrow \infty}\wtY^n$ exists.

\noindent {\it Step 1.} We will first prove that
\begin{align} \label{eq6.11}
\wtY_t \geq \essinf_{\tau\in\bSTFt}\esssup_{\sigma\in\STFt}\EE\big[\zeta_{\sigma}\I_\seq{\tau>\sigma}
+(\zeta\vee \eta)_\tau \I_\seq{\tau \leq \sigma}\,|\,\cF_t\big].
\end{align}
To establish \eqref{eq6.11}, for any fixed $t$ and $n$, we define $\bar{\tau}^n_t:= \inf\{s\in [t,T]:\wt{L}^n_s-\wt{L}^n_t > 0\}$ where $\wt{L}^n_t:=\int_{]0,t]} n(\wtY^n_s-\eta_s)^+\,dA_s$. Since the process $\wt{L}^n$ is continuous, the graph of the stopping time $\bar{\tau}^n_t$ is contained in $\barS \cap [t,T]$, which implies that $\bar{\tau}^n_t \in\bSTFt$. Suppose, on the contrary, that $\bar{\tau}^n_t \not \in\bSTFt$. Then the event $\{\bar\tau^n_t < T\}\cap \{\bar\tau^n_t \in \barS^c\}$ has a positive probability and, for any fixed $\omega\in \{\bar\tau^n_t < T\}\cap \{\bar\tau^n_t \in \barS^c\} $, there exists $\delta = \delta (\omega) > 0$ such that $A_{\bar\tau^n_t + \delta} =A_{\bar\tau^n_t}$. However, this contradicts the definition of $\bar\tau^n_t$ since $\wt{L}^n$ is absolutely continuous with respect to $\wt\Gamma$ and thus $\wt{L}^n_{\bar\tau^n_t + \delta} =\wt{L}^n_{\bar\tau^n_t}$.

The continuity of $A$ entails that $\wtY^n_{\bar{\tau}^n_t+} \geq \liminf_{s\downarrow \bar\tau^n_t}\eta_{s}$ on $\{\bar{\tau}^n_t<T\}$ and, consequently, using also  \eqref{eq6.9m} and the right-continuity of $\eta$ we deduce that $\wtY^n_{\bar{\tau}^n_t} = \wtY^n_{\bar{\tau}^n_t+}+
\Delta^+ \wt{K}^{n,g}_{\bar{\tau}^n_t}\geq \wtY^n_{\bar{\tau}^n_t+} \geq \eta_{\bar{\tau}^n_t}$.
In addition, we have $\wtY^n_{\bar{\tau}^n_t} \geq \zeta_{\bar{\tau}^n_t}$ since $(\wtY^n,\wtZ^n,\wtK^n)$ solves the reflected GBSDE \eqref{eq6.9m}.
We conclude that $\wtY^n_{\bar{\tau}^n_t} \geq (\zeta\vee \eta)_{\bar{\tau}^n_t}$ on $\{\bar{\tau}^n_t<T\}$ and, manifestly, $\wtY^n_T=\zeta_T$.

We now take an arbitrary stopping time $\sigma \in \STFt$ and define $\nu :=\bar{\tau}^n_t\wedge\sigma$ so that $\wtY^n$
is a strong supermartingale on $\llb t,\nu\rrb$ since $\wt{L}^n_{\nu}=\wt{L}^n_t$. Then $\wtY^n_\nu \geq \zeta_\nu$ on $E_1:=\{\bar{\tau}^n_t\geq \sigma\}$ and
$\wtY^n_{\nu}\geq (\zeta\vee \eta)_{\nu}$ on $E_2:=\{\bar{\tau}^n_t < \sigma\}$. Consequently, for any $\sigma\in\STFt$,
\begin{align}\label{eq6.11g}
\wtY^n_t\geq \EE[ \wtY^n_{\nu}\,|\,\cF_t]\geq \EE\big[ \zeta_\sigma\I_\seq{\bar{\tau}^n_t>  \sigma}
+(\zeta\vee \eta)_{\bar{\tau}^n_t}\I_\seq{\bar{\tau}^n_t \leq \sigma}\,|\,\cF_t\big]
\end{align}
from which we deduce that
\begin{align*}
\wtY^n_t & \geq \esssup_{\sigma\in\STFt} \EE\big[ \zeta_\sigma\I_\seq{\bar{\tau}^n_t>  \sigma}
+(\zeta\vee \eta)_{\bar{\tau}^n_t}\I_\seq{\bar{\tau}^n_t \leq \sigma}\,|\,\cF_t\big]\\
& \geq \essinf_{\tau\in\bSTFt}\esssup_{\sigma\in\STFt} \EE\big[ \zeta_\sigma\I_\seq{\tau>  \sigma}
+(\zeta\vee \eta)_{\tau}\I_\seq{\tau \leq \sigma}\,|\,\cF_t\big].
\end{align*}
Finally, the sequence $\wtY^n$ is decreasing and $\wtY=\lim_{\,n\rightarrow \infty}\wtY^n$ so that we obtain \eqref{eq6.11}.

\noindent {\it Step 2.}  In this step, we will establish the inequality
\begin{align} \label{eq6.12}
\wtY_t \leq \esssup_{\sigma\in\STFt}\essinf_{\tau\in\bSTFt}\EE\big[\zeta_{\sigma}\I_\seq{\tau > \sigma }
+(\zeta\vee \eta)_\tau\I_\seq{\tau \leq \sigma}\,|\,\cF_t\big]
\end{align}
by showing that, for any $\varepsilon > 0$, there exists $\bar{\sigma}_t\in\STFt$, which may depend on $\varepsilon$, such that for an arbitrary $\tau\in\bSTFt$ we have
\begin{align} \label{eq6.11b}
\wtY_t \leq \EE\big[\zeta_{\bar{\sigma}_t}\I_\seq{\tau > \bar{\sigma}_t}+(\zeta\vee \eta)_\tau \I_\seq{\tau \leq \bar{\sigma}_t}\,|\,\cF_t\big] + \varepsilon.
\end{align}

For a fixed $t$ and $\varepsilon>0$, we define $\bar\sigma^n_t := \inf \{ s\in [t,T]: \wtY^n_s \leq \zeta_s + \varepsilon\}$. Recall that the sequence $\wtY^n$ is monotonically decreasing as $n\rightarrow \infty$ and $\wtY=\lim_{\,n\rightarrow \infty}\wtY^n$ so that $\bar\sigma^n_t \geq \bar\sigma^{n+1}_t$. We define an $\FF$-stopping time $\bar\sigma_t:= \lim_{\,n\rightarrow \infty} \bar\sigma^n_t$.  From the lower bound in \eqref{eq6.11g} we know that $\wtY^n_t \geq 0$, while the comparison theorem for reflected GBSDEs gives, for every $n\in\NN$,
\begin{align} \label{eq6.11x}
\wtY_t \leq \wtY^n_t\leq X_t=\zeta_T-\int_{]t,T]} Z_s\,dM_s+K_T-K_t=\esssup_{\tau \in \STFt} \EE[\zeta_\tau\,|\,\cF_t] \leq c_\zeta
\end{align}
where $(X,Z,K)$ is a solution to the reflected BSDE implicit in \eqref{eq6.11x} with the lower obstacle $\zeta$ (see Section \ref{nsec3.1}) and thus the second equality is due to the well-known relationship between a solution to the reflected BSDE with null generator and the value process of an optimal stopping problem with the bounded reward process $\zeta$.

From the assumption that $\zeta$ is right-upper-semicontinuous we deduce that $\wtY^n_{\bar\sigma^n_t} \leq \zeta_{\bar\sigma^n_t} + \varepsilon$ where the inequality is trivially satisfied on the event $\{\bar\sigma^n_t = T\}$. Since $\wtY^n$ is a l\`adl\`ag process, it would be possible to have sample paths satisfying: $\wtY^n_{\bar{\sigma}^n_t} > \zeta_{\bar\sigma^n_t}+\varepsilon$ and there exists $\delta >0$ such that $\wtY^n \leq \zeta+\varepsilon$ on $\rrb\bar\sigma^n_t , \bar\sigma^n_t +\delta\rrb$. However, by the right-upper-semicontinuity of $\zeta$, this would imply the inequalities $\wtY^n_{\bar\sigma^n_t+} \leq \zeta_{\bar\sigma^n_t}+\varepsilon$ and $\Delta^+ \wtK^{n,g}_{\bar\sigma^n_t}>0$. This would lead to a contradiction since, from the Skorokhod condition for $\wtK^{n,g}$, the inequality $\Delta^+ \wtK^{n,g}_{\bar\sigma^n_t}>0$ implies that $\wtY^n_{\bar\sigma^n_t} = \zeta_{\bar\sigma^n_t} < \zeta_{\bar\sigma^n_t} + \varepsilon < \wtY^n_{\bar\sigma^n_t}$. In view of these considerations, we conclude that $\wtY^n > \wtY^n -\varepsilon> \zeta$ on $\llb t,\bar\sigma^n_t \llb$ and $\wtY^n_- > \wtY^n_- -\varepsilon\geq  \bar \zeta$ on $\rrb t,\bar\sigma^n_t\rrb$. Together with the Skorokhod condition for the process $\wtK^{n}$, this gives
\begin{align} \label{kzero}
\wtK^n_{\bar\sigma^n_t}-\wtK^n_t=\int_{\rrb t,\bar\sigma^n_t\rrb}d\wtK^{n,r}_s + \int_{\llb t,\bar\sigma^n_t\llb} d\wtK^{n,g}_{s+} = 0.
\end{align}
If we take $\nu:=\tau \wedge \bar{\sigma}^n_t$ where $\tau\in\bSTFt$ is arbitrary, then
\begin{align*}
\wtY^n_t& = \EE\big[\wtY^n_{\nu}-\int_{\rrb t,\nu\rrb} n(\wtY^n_s-\eta_s)^+\,dA_s+\wt{K}^n_{\nu}-\wt{K}^n_t\,\big|\,\cF_t\big]\leq \EE\big[\wtY^n_{\nu}|\,\cF_t\big] \nonumber\\
& = \EE\big[\wtY^n_{\bar{\sigma}^n_t}\I_\seq{\tau > \bar{\sigma}^n_t}
          +\wtY^n_{\tau}\I_\seq{\tau \leq  \bar{\sigma}^n_t} \,|\,\cF_t\big] \nonumber\\
& \leq \EE\big[\zeta_{\bar{\sigma}^n_t}\I_\seq{\tau > \bar{\sigma}^n_t}+
     (\wtY^n\vee \zeta\vee \eta)_{\tau}\I_\seq{\tau \leq \bar{\sigma}^n_t}\,|\,\cF_t\big] + \varepsilon\\
& \leq      \EE\big[\zeta_{\bar{\sigma}^n_t}\I_\seq{\tau > \bar{\sigma}^n_t}+
     (Y\vee \zeta\vee \eta)_{\tau}\I_\seq{\tau \leq \bar{\sigma}_t}\,|\,\cF_t\big] + \EE\big[|\wtY^n_\tau-Y_\tau|\,|\,\cF_t\big] + C\EE\big[\I_{\rrb\bar\sigma_t, \bar\sigma^n_t \rrb}(\tau)|\,\cF_t\big] + \varepsilon
\end{align*}
where on the event $E_1:= \{\tau > \bar\sigma^n_t\}$ we have used the inequality $\wtY^n_{\bar{\sigma}^n_t} \leq \zeta_{\bar{\sigma}^n_t} + \varepsilon$ while on the event $E_2:= \{\tau \leq \bar\sigma^n_t\}$ we have used the trivial inequality $\wtY^n_{\tau}\leq\wtY^n_{\tau}\vee \zeta_{\tau}\vee \eta_{\tau}$.
By considering the limit superior in $n$ and using the conditional reverse Fatou lemma (see Theorem 2 in \cite{Z1998})  together with the right-upper-semicontinuity of $\zeta$ and the monotone convergence theorem, we obtain
\begin{align*}
\wtY_t \leq 	\EE\big[\zeta_{\bar{\sigma}_t}\I_\seq{\tau > \bar{\sigma}_t}+
     (Y\vee \zeta\vee \eta)_{\tau}\I_\seq{\tau \leq \bar{\sigma}_t}\,|\,\cF_t\big] +\varepsilon.
\end{align*}
Next, we will show that $Y$ can be omitted from the conditional expectation above. For any $\tau \in \bSTFt$, equation \eqref{kzero} gives
\begin{align*}
\wtY^n_{\tau}=  \wtY^n_{\bar\sigma^n_\tau}-\int_{\rrb\tau,\bar\sigma^n_\tau\rrb} \wtZ^n_s\,dM_s-\int_{\rrb \tau, \bar\sigma^n_\tau\rrb} n(\wtY^n_s-\eta_s)^+\,dA_s.
\end{align*}

We now use similar arguments as in Step 2 in the proof of Theorem \ref{th6.1}. We observe that $\wtY^n_{\bar\sigma^n_\tau}\leq \zeta_{\bar\sigma^n_\tau}+\varepsilon$ and, for all $(\omega,s, y) \in \Omega\times [0,T]\times \RR$,
$$
(y-\eta_s)^+(\omega) \geq \I_{\llb 0, \bar\sigma_\tau\rrb}(s) (y-\eta_s)^+(\omega) \geq \I_{\llb 0, \bar\sigma_\tau\rrb}(s) (y-\eta_s)(\omega)
$$
where the function $g(\omega,s, y):=\I_{\llb 0, \bar\sigma_\tau\rrb}(s) (\eta_s-y)(\omega)$ is nonincreasing in $y$, for every $(\omega,s)\in \Omega\times [0,T]$.  By applying the comparison theorem for a GBSDE (see Proposition \ref{npro1.1}) on the interval $\llb \tau, \bar\sigma^n_\tau\rrb$, we see that $\wtY^n \leq Y^n$ where $(Y^n,Z^n)$ solves the following linear BSDE
\begin{align*}
			 Y_\tau^n  & =  \zeta_{\bar\sigma^n_\tau}  + \varepsilon-\int_{\rrb\tau,\bar\sigma^n_\tau\rrb} Z^n_s\,dM_s-\int_{\rrb \tau, \bar\sigma^n_\tau\rrb} n(Y^n_s-\eta_s)\I_{\llb 0, \bar\sigma_\tau\rrb}(s) \,dA_s\\
			 & =  \zeta_{\bar\sigma^n_\tau}  + \varepsilon-\int_{\rrb\tau,\bar\sigma^n_\tau\rrb} Z^n_s\,dM_s+ \int_{\rrb \tau, \bar\sigma^n_\tau\rrb} n(\eta_s- Y^n_s)\,dA^{\bar\sigma_\tau}_s.
\end{align*}
Let us denote $A^n := nA$. Since $\tau \leq \bar\sigma_\tau \leq \bar\sigma^n_\tau \leq T$, Corollary \ref{ncor1.2} gives
\begin{align*}
 Y_\tau^n& = \EE\big[(\zeta_{\bar\sigma^n_\tau}+ \varepsilon)\cE_{\tau,\bar\sigma_\tau}(-A^n)+(\I_{\rrb \tau,\bar\sigma_\tau\rrb}\eta\cE_{\tau,\cdot}(-A^n)\bcd A^n)_{T}\, | \cF_\tau\big]\\
& \leq \EE\big[\zeta_{\bar\sigma_\tau}\cE_{\tau,\bar\sigma_\tau}(-A^n)+(\I_{\rrb \tau,\bar\sigma_\tau\rrb}\eta\cE_{\tau,\cdot}(-A^n)\bcd A^n)_{T}\, | \cF_\tau\big]  + \EE\big[[\zeta_{\bar\sigma^n_\tau} - \zeta_{\bar\sigma_\tau}]\cE_{\tau,\bar\sigma_\tau}(-A^n)  |\cF_\tau\big] + \varepsilon
\end{align*}
where in the last inequality we have used the inequality $\cE_{\tau,\bar\sigma_\tau}(-A^n) \leq 1$. The quantity $\cE_{\tau,\bar\sigma_\tau}(-A^n)$ converges to $\I_{\{\tau=\bar\sigma_\tau\}}$ as $n\rightarrow \infty$ and, by the subadditivity of the limit superior, the conditional reverse Fatou lemma and the dominated convergence theorem, we obtain
\begin{align*}
&\limsup_{\,n\rightarrow \infty}\, \EE\big[[\zeta_{\bar\sigma^n_\tau} - \zeta_{\bar\sigma_\tau}]\cE_{\tau,\bar\sigma_\tau}(-A^n) \,|\,\cF_\tau\big]\\
& \leq \limsup_{\,n\rightarrow \infty}\, \EE\big[[\zeta_{\bar\sigma^n_\tau} \cE_{\tau,\bar\sigma_\tau}(-A^n)\,|\,\cF_\tau\big] - \lim_{\,n\rightarrow \infty}\, \EE\big[\zeta_{\bar\sigma_\tau}\cE_{\tau,\bar\sigma_\tau}(-A^n)\,|\,\cF_\tau\big]\\
& \leq  \EE\big[[\limsup_{\,n\rightarrow \infty} \zeta_{\bar\sigma^n_\tau} - \zeta_{\bar\sigma_\tau}]\I_{\{\tau=\bar\sigma^n_\tau\}}\,|\,\cF_\tau\big] \leq 0
\end{align*}
where the last inequality holds since $\zeta$ is right-upper-semicontinuous along stopping times (see Remark B.3 in \cite{KQC2014}). For any fixed $\varepsilon > 0$, we conclude from the subadditivity of the limit superior and Lemma \ref{lem3.4x} that
\begin{align*}
\wtY_\tau\leq \zeta_{\tau}\I_{\{\tau=\bar\sigma_\tau\}}+\eta_{\tau}\I_{\{\tau<\bar \sigma_\tau\}}+\varepsilon \leq (\zeta\vee \eta)_\tau + \varepsilon
\end{align*}
and thus $\wtY_\tau \leq (\zeta\vee \eta)_\tau$ for every stopping time $\tau$ in $\bSTFt$, which gives the desired upper bound in \eqref{eq6.11b}.

\noindent {\it Step 3.} Since we always have that
\begin{align*}
\essinf_{\tau\in\bSTFt}\esssup_{\sigma\in\STFt}\EE[\Theta(\sigma,\tau)\,|\,\cF_t]\geq \esssup_{\sigma\in\STFt} \essinf_{\tau\in\bSTFt}\EE[\Theta(\sigma,\tau)\,|\,\cF_t]
\end{align*}
we obtain \eqref{eq6.7} by combining \eqref{eq6.11} with \eqref{eq6.12}.
\end{proof}

\brem  \label{remx.x}
It is natural to conjecture that the process $\wtY$ given by \eqref{eq6.7} can be represented through a solution to a doubly reflected BSDE.
Although we do not examine that issue in the present work, let us point out that in the case where the processes $\eta$ and $\zeta$ are c\`adl\`ag and $A_t = t$, it was demonstrated in Theorem 3.1 of Hamad\`ene et al. \cite{HHO2010} that the limit $\wtY$ satisfies locally a doubly reflected BSDE
with the lower and upper obstacle equal to $\zeta$ and $\eta$, respectively. In addition, Theorem 4.1 in \cite{HHO2010} shows that if the obstacles are completely separated, in the sense that the inequality $\eta > \zeta$ holds, then $\wtY$ is a solution to a doubly reflected BSDE.

From the point of view of the present work, we expect analogous results to be valid. However, since the exercise set of the minimizer in \eqref{eq6.7}  is constrained to the right support of the measure generated by $A$, we expect the limit $\wtY$ to satisfy a doubly reflected BSDE with possible constraints on the Skorokhod condition. More specifically, we expect that the reflection process associated with the upper obstacle $\eta$ only increases at times when a solution $\wtY$ hits the upper obstacle $\eta$ and, simultaneously, the process $A$ increases (more precisely, on the right support of the random measure generated by $A$).
\erem

%%%%%%%%%%%%%%%%%%%%%%%%%%%%%%%%%%%%%%%%%%%%%%%%%%%%%%%%%%%%%%%%%%
%%%%%%%%%%%%%%%%%%%%%%%%%%%%%%%%%%%%%%%%%%%%%%%%%%%%%%%%%%%%%%%%%%
%%%%%%%%%%%%%%%%%%%%%%%%%%%%%%%%%%%%%%%%%%%%%%%%%%%%%%%%%%%%%%%%%%

\end{document}